\documentclass[10pt, sigconf]{acmart}
\setcopyright{none}
\usepackage{etoolbox}
\makeatletter
\patchcmd{\maketitle}{\@copyrightspace}{}{}{}
\makeatother
\renewcommand\footnotetextcopyrightpermission[1]{} 

\usepackage{graphicx}
\usepackage{dsfont} 
\usepackage{amsmath}
\usepackage{amssymb}
\usepackage{url}
\usepackage{amsthm}

\usepackage{multicol}

\usepackage{color}
\usepackage{subcaption}

\newcommand{\doi}[1]{\href{http://dx.doi.org/#1}{\normalsize{\textsc{doi:}}~\nolinkurl{#1}}}
\newcommand{\arxiv}[1]{\href{http://arxiv.org/abs/#1}{\normalsize{\textsc{arxiv:}}~\nolinkurl{#1}}}

\renewcommand{\epsilon}{\varepsilon}
\renewcommand{\phi}{\varphi}

\newcommand{\R}{\mathbb{R}}

\newcommand{\N}{\mathbb{N}}

\newcommand{\ubar}[1]{\text{\b{$#1$}}}

\let\originalleft\left
\let\originalright\right
\renewcommand{\left}{\mathopen{}\mathclose\bgroup\originalleft}
\renewcommand{\right}{\aftergroup\egroup\originalright}

\def\clap#1{\hbox to 0pt{\hss#1\hss}}

\newcommand{\norm}[1]{\left\lVert #1\right\rVert}

\newcommand{\set}[1]{\left\{ #1\right\}}

\let\p=\paren

\newcommand{\sqparen}[1]{\left[ #1\right]}

\let\sp=\sqparen

\newcommand{\derivoper}{\mathrm{d}}

\newcommand{\deriv}[3][]{%
  \ifx\relax#1\relax{
    \frac{\derivoper #2}{\derivoper #3}%
  }\else{%
    \frac{\derivoper^{#1} #2}{\derivoper #3^{#1}}%
  }\fi%
}
\newcommand{\pderiv}[3][]{%
  \ifx\relax#1\relax{
    \frac{\partial #2}{\partial #3}%
  }\else{%
    \frac{\partial^{#1} #2}{\partial #3^{#1}}%
  }\fi%
}
\newcommand{\grad}[2][]{%
  \ifx\relax#1\relax{
    \nabla #2%
  }\else{%
    \nabla_{#1} #2%
  }\fi%
}
\newcommand{\diver}[2][]{%
  \ifx\relax#1\relax{
    \nabla \cdot #2%
  }\else{%
    \nabla_{#1} \cdot #2%
  }\fi%
}
\newcommand{\curl}[2][]{%
  \ifx\relax#1\relax{
    \nabla \times #2%
  }\else{%
    \nabla_{#1} \times #2%
  }\fi%
}
\newcommand{\lapl}[2][]{%
  \ifx\relax#1\relax{
    \Delta #2%
  }\else{%
    \Delta_{#1} #2%
  }\fi%
}

\newcommand{\smats}[1]{%
  \sqparen{\begin{smallmatrix}#1 \end{smallmatrix}}%
}

\usepackage{enumitem}

\setlist[enumerate]{leftmargin=*}
\setlist[itemize]{leftmargin=*}

\newcommand{\indmetric}[1]{\widetilde{d}_{#1}}

\usepackage{cleveref}
\newcommand{\hs}{\mathcal{H}}
\newcommand{\hse}{\mathcal{H}^\epsilon}

\newcommand{\M}{\mathcal{M}}
\newcommand{\Me}{\mathcal{M}^\epsilon}
\newcommand{\J}{\mathcal{J}}
\newcommand{\BVU}{PC\p{\sp{0,T}, U}}
\theoremstyle{acmplain}
\newtheorem{assumption}{Assumption}

\makeatletter
\def\@copyrightspace{\relax}
\makeatother

\begin{document}

\title{%
  On the Relaxation of Hybrid Dynamical Systems %
}

\author{ Tyler Westenbroek \and S. Shankar Sastry \and Humberto Gonzalez }

\begin{abstract}
Hybrid dynamical systems have proven to be a powerful modeling abstraction, yet fundamental questions regarding the dynamical properties of these systems remain. In this paper, we develop a novel class of relaxations which we use to recover a number of classic systems theoretic properties for hybrid systems, such as existence and uniqueness of trajectories, even past the point of Zeno. Our relaxations also naturally give rise to a class of provably convergent numerical approximations, capable of simulating through Zeno. Using our methods, we are also able to perform sensitivity analysis about nominal trajectories undergoing a discrete transition -- a technique with many practical applications, such as assessing the stability of periodic orbits. 
\end{abstract}

\maketitle

\section{Introduction}

While hybrid dynamical systems have proven to be a highly expressive modeling framework, the flexibility they provide does not come without its challenges. Despite considerable efforts to extract classic systems theoretic properties from hybrid systems in works such as \cite{simic2000towards} and \cite{lygeros2003dynamical}, fundamental questions regarding even the existence and uniqueness of their executions abound, as the interplay between their discrete and continuous dynamics is not fully understood. 

Perhaps the most notable phenomena unique to hybrid systems, Zeno executions \cite{RNC:RNC592} arise when an infinite number of discrete transitions occur in a finite amount of time. In order to accommodate Zeno trajectories into theoretical and computational frameworks, a number of techniques have been proposed. In \cite{johansson1999regularization}, the authors propose techniques to regularize hybrid systems in time or space, which prevent an infinite number of transitions from occurring. Yet they are able to prove convergence for their relaxations only for Zeno executions which accumulate to a single point. In \cite{burden2015metrization}, the authors extend this proof of convergence to the numerical setting. Alternatively, the authors in \cite{or2009formal} go to great lengths to identify Zeno executions, and replace them with executions of a reduced order dynamical system, in order to avoid directly handling Zeno. The authors are able to extend simulations past the point of Zeno in some cases, but the results only hold for mechanical systems.

Even if we disregard the pathologies introduced by Zeno executions, a number of theoretical and practical challenges remain to fully understand the executions of hybrid systems. The trajectories of hybrid systems may be discontinuous with respect to inputs and initial conditions \cite{lygeros2003dynamical}, and in such cases may not be faithfully approximated in the numerical setting. Indeed, many works focussed on numerical integration for hybrid systems such as \cite{burden2015metrization} and \cite{Esposito2001} make restrictive assumptions about the trajectories being simulated to accommodate this obstacle, and require that timesteps be placed in small neighborhoods around discrete events. 

Taking steps to overcome these limitations, we introduce a novel relaxation scheme for hybrid dynamical systems. First, we demonstrate how to reduce a discrete jump of a hybrid system to the execution of a switched system. This enables us to use the solution concept of Filippov \cite{filippov2013differential} to define closed form solutions for some Zeno executions of our hybrid systems. We then extend the procedure presented in~ \cite{llibre2007regularization} to regularize this collection of switched systems, which recover the sliding solutions of Filippov in the limit, and use the resulting vector fields to construct trajectories over our \emph{relaxed hybrid systems}. To construct these relaxations we take the approach of adding an epsilon-thick strip to each of our guard sets as in \cite{johansson1999regularization}, and endow our relaxed hybrid systems with the topology from \cite{burden2015metrization}.

Using this framework we extend the state of the art in several directions. Firstly, we use the limit of our relaxations to construct a novel solution concept for hybrid systems, wherein the trajectory generated by ever pair of initial conditions and inputs is unique and well defined, even past the point of Zeno. Secondly, our relaxation yields a provably convergent numerical approximations, which can approximate solutions past Zeno. Finally, we are able to perform sensitivity analysis on the trajectories of our relaxations -- which has numerous practical applications, which we discuss in our closing remarks.


\section{Filippov Solutions}

\begin{figure}[t!]
  \includegraphics[width= .5\columnwidth, height = .95\textheight, keepaspectratio ]{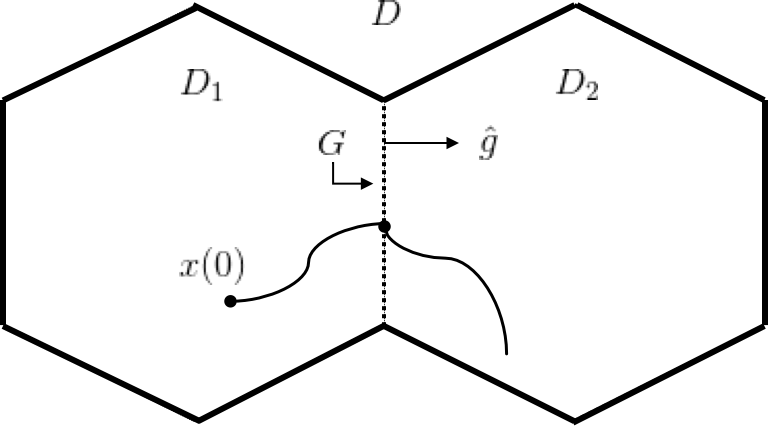}%
  \caption{A bimodal switched system with a representative trajectory $x$.}
  \label{fig:switched}
\end{figure}

We briefly introduce the solution concept of Filippov \cite{filippov2013differential} for switched systems. Consider the bimodal switched system depicted in Figure~\ref{fig:switched}, where the domains $D_1$, $D_2 \subset \R^n$ are separated by the plane $G = \set{ x \in \R^n \colon g(x) = \hat{g}^T x - c =0 }$, where $\hat{g} \in \R^n$ is a unit vector and $c$ is a scalar. Defining $D = D_1 \cup D_1$, and an allowable set of inputs $U$, the dynamics of the system are governed by $f \colon D \times U \to \R^n$ where
\begin{equation}
f(x,u) = \begin{cases}
f_1(x,u), & \text{if } x \in D_1 \\
f_2(x,u), & \text{if } x \in D_2 ,\\
\end{cases}
\end{equation}
and $f_1\colon D_1 \times U \to \R^n$ and $f_2 \colon D_2 \times U \to \R^n$ are both continuous but  $f$ may be discontinuous along $G$. In particular, we align $\hat{g}$ such that $g_e(x) \leq 0$, $\forall x \in D_1$, and $g_e(x) \geq 0$, $\forall x \in D_2$.

We may partition the set $H = G \times U$ into three distinct regions: the \emph{crossing region} ($\Sigma_c$), the \emph{sliding region} ($\Sigma_s$), and the \emph{escaping region} ($\Sigma_{es}$). We characterize these regions by $\Sigma_c  =\set{(x,u) \in H \colon \p{\hat{g}^Tf_1(x,u)}\p{\hat{g}^Tf_2(x,u)} > 0}$, $\Sigma_s  =\set{(x,u) \in H \colon \p{\hat{g}^Tf_1(x,u)} > 0, \p{\hat{g}^Tf_2(x,u)} < 0} $, and

\noindent
$\Sigma_{es}=\set{(x,u) \in H \colon \hat{g}^Tf_1(x,u)\leq 0, \hat{g}^Tf_2(x,u)  \geq 0}$.
  
When $(x,u) \in \Sigma_c$, the trajectory $x$ simply crosses from one domain to the other. When $(x,u) \in \Sigma_s$, both $f_1$ and $f_2$ are pointing into the surface $G$, confining the trajectory to this set. One way to model trajectories in this regime is to switch between the vector fields $f_1$ and $f_2$ infinitely fast -- i.e as a Zeno execution. However, Filippov solutions offer us another route to understand such systems. The \emph{Filippov sliding vector field} for this system, $f^s \colon \Sigma_s \to \R^n$, is given by $f^s(x,u) = \p{1- \alpha(x,u)} f_1(x,u) + \alpha(x,u) f_2(x,u) $,
where $\alpha \colon \Sigma_s \to \sp{0,1}$ is defined by
$\alpha(x,u) = \frac{\hat{g}^T f_1(x,u) }{\hat{g}^T f_1(x,u) - \hat{g}^T f_2(x,u)}$.
In particular, $\alpha$ is constructed such that $\forall (x,u) \in \Sigma_s$, we have $\hat{g}^Tf^s(x,u) =0$, confining the solution to $G$ as desired. Thus, we may use the solution concept of Filippov replace some Zeno trajectories with well defined vector fields. However, when $(x,u) \in \Sigma_{es}$, the solution concept of Filippov breaks down, and there are Zeno trajectories that are left ill-defined. 

\section{Hybrid Systems}
We introduce our class of hybrid systems, inspired by \cite{burden2015metrization}.
\begin{definition} \label{def:hybrid_sys}
A hybrid system is a seven-tuple \\  $\mathcal{H} =  \p{\mathcal{J}, \Gamma, \mathcal{D}, \mathcal{U}, \mathcal{F}, \mathcal{G}, \mathcal{R}}$, where:
\begin{itemize}
\item $\mathcal{J}$ is a finite set indexing the discrete states of $\mathcal{H}$; 
\item $\Gamma \subset \mathcal{J} \times \mathcal{J}$ is the set of edges, forming a graphical structure over $\mathcal{J}$, where edge $e = (j,j') \in \Gamma$ corresponds to a transition from $j$ to $j'$;
\item $\mathcal{D} =  \set{D_j}_{j \in \mathcal{J}}$ is the set of domains, where $D_j$ is a compact $n$-dimensional polytope in $\R^{n}$, $n \in \N$; \label{item:domain_def}
\item $U \subset \R^{m}$ is a compact set of inputs, $m \in \N$;
\item $\mathcal{F} = \set{f_j}_{j \in \mathcal{J}}$ is the set of vector fields, where each $f_j \colon D_j \times U \to \R^{n}$ is continuously differentiable \footnote{This ensures continuous state trajectories are unique and well defined on our continuous domains, since continuous functions are Lipschitz over compact sets.} and defines the continuous dynamics of the system on $D_j$\footnote{We incur no loss of generality by considering time invariant vector fields. Indeed, one may add time as a continuous state $z$, with dynamics $\dot{z} =1$ and initial condition $z(t_0) =t_0$.}  ; 
\item $\mathcal{G} = \set{G_e}_{e =(j,j') \in \Gamma}$ is the set  of guards, where each $G_{(j,j')} \subset \partial D_{j}$ is a codimension 1 plane with corner; that is, there exists a unit vector $\hat{g}_e \in \R^n$ and a scalar $c_e$ such that $G_e \subset \set{x \in \partial D_j \colon g_e(x) = \hat{g}_e^Tx-c_e = 0}$\footnote{We choose the convention that $\hat{g}_e$ 'points out' of $D_j$ along $G_e$ -- i.e. $g_e(x) \leq 0, \ \forall x \in D_j$. }; and,\label{linear_guards}
\item $\mathcal{R} = \set{R_e}_{e = (j,j') \in \Gamma}$ is the set of reset maps where, for each $e \in \Gamma$, $R_{e} \colon G_{e} \to \partial  D_{j'}$ is defined by $R_e(x) = A_ex + b_e$, where $A_e \in \R^{n \times n}$ and $b_e \in \R^n$.

\end{itemize}
\end{definition} 
When the guard $G_{(j,j')}$ is crossed, a discrete transition from mode $j$ to $j'$ occurs, and the continuous state is instantaneously rest by $R_{(j,j')}$. We unify our continuous and discrete state spaces using the concept of a \emph{disjoint union}. That is, we embed our continuous domains in the space $\coprod_{j \in J}D_j = \bigcup_{j \in \J} D_j \times \set{j}$. By an abuse of notation, throughout the paper we shall simply use $D_j$ to refer to $D_j \times \set{j}$. For each $j \in \J$, we let $\mathcal{N}_j = \set{e \in \Gamma \colon \exists j' \in \J \ s.t. \ e = (j,j')}$ be the \emph{neighborhood} of $\J$. We let $\BVU$ denote the class of piecewise continuous functions from the interval $\sp{0,T}$ to $U$.  For each $j \in \J$, let $f_j \colon \R^n \times U \to \R^n$ be any continuously differentiable extension to $f_j \colon D_j \times U \to \R^n$, guaranteed to exist by Lemma 5.6 of \cite{lee2003smooth}. Abusing notation, throughout the paper, when the symbol $f_j$ is used, it is understood that we are referring to the extended version of the function. Similarly, for each $e \in \Gamma$ we extend $R_e \colon \R^n \to \R^n$, where for each $x \in \R^n$ we still define $R_e(x) = A_ex +b_e$. We impose the following assumptions to simplify the discussion sliding vector fields throughout the paper.

\begin{assumption}\label{ass:reversible}
Let $e \in \Gamma$. Then $A_e$ is invertible, and there exists an edge $e' \in \Gamma$ such that $R_e(G_e) = G_{e'}$, $R_{e'}(G_{e'}) = G_e$, and, $\forall x \in G_e$, $x = R_{e'}(R_e(x))$. 
\end{assumption}

We say that the edge $e$ is \emph{reversible} if it satisfies  Assumption~\ref{ass:reversible}, since a transition along $e$ can be 'reversed' by a transition along $e'$, which we refer to as the \emph{partner} of $e$. We will often use $e'$ to refer to the partner of $e$ with out explicitly stating their relationship. The one-to-one correspondence between $G_{e}$ and $G_{e'}$ defined by $R_e$ and $R_{e'}$ will simplify our initial discussion of Filippov solutions along the guard sets of our hybrid systems. In the optional appendix, we outline how to overcome this assumption in theory, and in Section~\ref{sec:examples} we produce numerical examples where the edges are not reversible. However, in both cases we only consider Zeno trajectories involving at most two edges of a hybrid system. 
\begin{assumption}\label{ass:overlapping}
For each pair of edges $e \neq \bar{e}$, $G_{e} \cap G_{\bar{e}} =\emptyset$.
\end{assumption}
While much work has been done to extend the solution concept of Filippov to cases where multiple continuous domains interface \cite{dieci2011sliding}, many open problems regarding the existence and uniqueness of solutions in such cases remain, and we wish to avoid such questions here in favor of presenting the main conceptual and technical components of our relaxation framework. We are currently investigating ways to extend our results to hybrid systems with non-linearities in their guard sets and reset maps, and overlapping guards. We now endow our hybrid systems with the topology from \cite{burden2015metrization}, which uses the concept of a \emph{quotient space} [\cite{kelley2017general}, Ch. 3]. 

Given a topological space $\mathcal{S}$ and a function $f \colon A \to B$,  where $A,B \subset \mathcal{S}$, we define the following equivalence relation:  $A \sim B = \set{(a,b) \in \mathcal{S}\times \mathcal{S} \colon a \in f^{-1}(b)}$,
and denote the quotient of $\mathcal{S}$ under $A \sim B$ by $\frac{\mathcal{S}}{\Lambda_f}$. Quotienting a space is often informally referred to as applying "topological glue" -- that is, for each $a \in A$ the sets $a$ and $f(a)$ are "glued" together, becoming a single set in $\frac{\mathcal{S}}{\Lambda_f}$. We embed our hybrid systems in the quotient space defined by their collection of reset maps. For hybrid system $\hs$ define
$\widehat{R}\colon \coprod_{e\in \Gamma} G_e \to \coprod_{j  \in \J} \partial D_{j'}$
by $\widehat{R}(x) = R_e(x)$ for each $x \in G_e$. Then the hybrid quotient space of $\hs$ is $
\M = \frac{\coprod_{j \in J} D_j }{\Lambda_{\widehat{R}}}$.

\begin{figure}[t]
  \centering
  \includegraphics[width= .7\columnwidth, height = .95\textheight, keepaspectratio ]{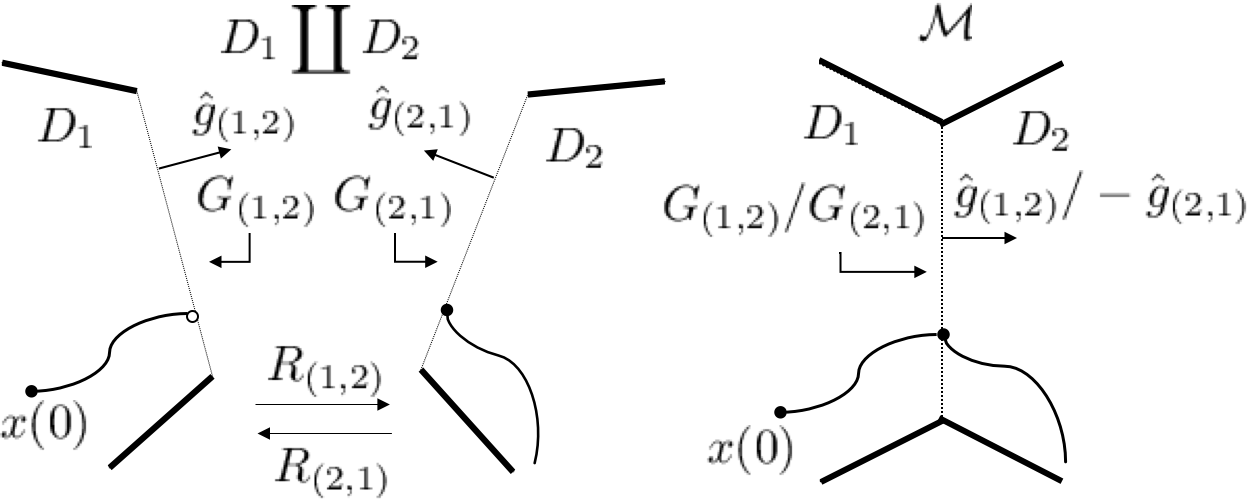}%
  \caption{A hybrid trajectory $x$ transitions from $D_1$ to $D_2$, represented on both $D_1 \coprod D_2$ (left) and $\M$ (right).}
  \label{fig:hybrid_quotient_space}
\end{figure}

The construction of $\M$ for a bimodal hybrid system $\hs$ is depicted in Figure \ref{fig:hybrid_quotient_space}, wherein the trajectory $x$ undergoes a discrete transition. Note, that for partners (1,2) and (2,1), the sets $G_{(1,2)}$ and $G_{(2,1)}$, while disjoint in $D_{(1,2)} \coprod D_{(2,1)}$, compose a single \emph{hybrid surface} in $\M$. Thus, the trajectories of our hybrid systems, to be defined in Section~\ref{sec:executions}, are in fact continuous on this space~\cite{burden2015metrization}.
 Speaking informally, the construction of $\M$ reduces $\hs$ into a switched system where the single hybrid surface $G_{(1,2)} \slash G_{(2,1)}$ separates $D_1$ and $D_2$. In section \ref{sec:switched_systems}, we do in fact demonstrate how to represent a discrete transition on $\M$ using the execution of a switched system. This will empower us to use the solution concept of Filippov to construct Zeno (sliding) trajectories along the hybrid surface $G_{(1,2)} \slash G_{(2,1)}$.

To understand the main difficulty in accomplishing this task, note that, when we construct $\M$ and describe a transition along the edge $(1,2)$, we apply an implicit change of coordinates wherein $i)$ we align the vector $\hat{g}_{(1,2)}$ with the vector $-\hat{g}_{(2,1)}$, so as a trajectory leaves $D_1$ it flows to the interior of $D_2$, and $ii)$ the matrix $A_{(1,2)}$ defines a correspondence between the surfaces $G_{(1,2)}$ and $G_{(2,1)}$; that is, $A_{(1,2)}$ transforms vectors in the $(n-1)$ dimensional subspace parallel to $G_{(1,2)}$ to lie in the $(n-1)$-dimensional subspace parallel to $G_{(2,1)}$. In Section~\ref{sec:switched_systems} we make this transformation explicit. 

Formally, in order metricize the hybrid quotient space we employ the induced length metrics from \cite{burden2011numerical}. Let $d \colon \R^n \times \R^n \to \R$ be a metric, then define $\tilde{d}_{\M} \colon \M \times \M \to \R$ for each $x, y \in \M$ by
\begin{equation}\label{eq:metric}
\tilde{d}_{\M}(x,y) = \inf_{k \in \N}\set{\sum_{i = 1}^k d(p_i,q_i) \colon x = p_1, \ y = q_k, \ q_i \sim p_{i+1}}.
\end{equation}
Intuitively, $\tilde{d}_{\M}(x,y)$ is simply the length of the shortest curve between $x$ and $y$ on $\M$, which may traverse multiple edges to connect the two points. 
 
\section{Relaxed Hybrid Systems}\label{def:relaxed_hybrid_system}
\begin{figure}[t]
  \centering
  \includegraphics[width= .9\columnwidth, height = .95\textheight, keepaspectratio ]{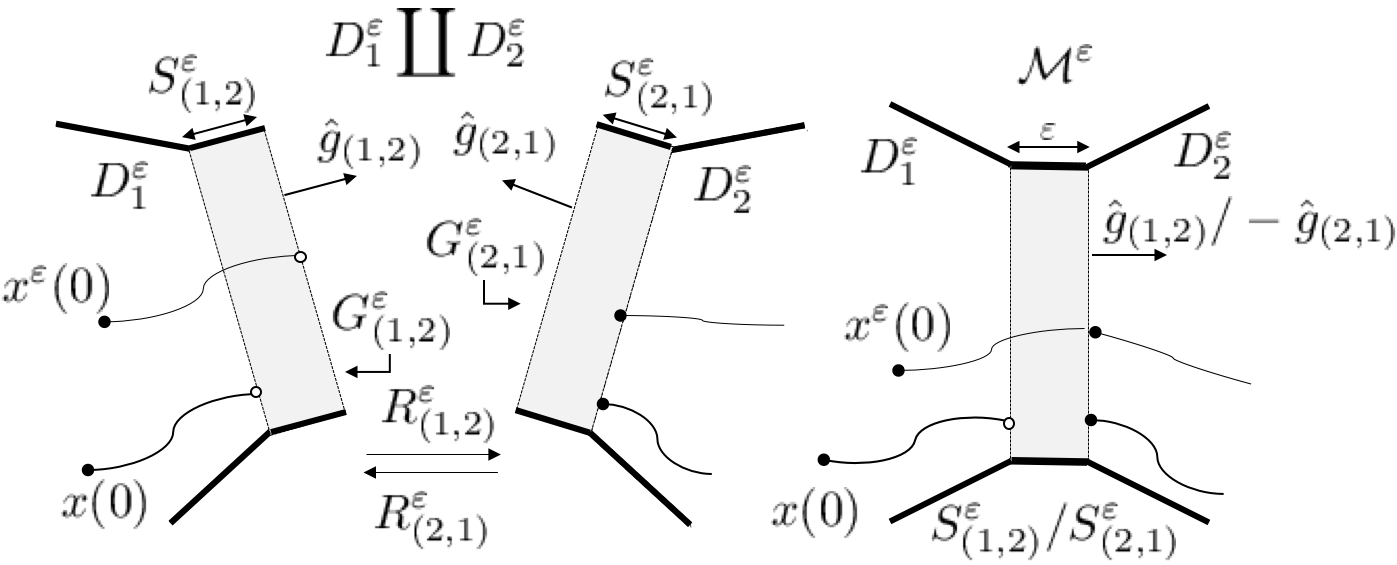}%
  \caption{A hybrid trajectory $x$ and relaxed hybrid strajectory $x^\epsilon$ transition from mode 1 to mode 2, represented on both $D_1^\epsilon \coprod D_2^\epsilon$ (left) and $\Me$ (right).}
  \label{fig:relaxed_hybrid_quotient_space}
\end{figure}

We now produce our definition for \emph{relaxed hybrid systems} inspired by \cite{burden2015metrization}, attaching an $\epsilon$-thick strip to each of the guard sets. In Section~\ref{sec:regularized_vfields}, we will define the continuously differentiable vector fields that we impart over these strips, which we will use to approximate Zeno trajectories. The relaxation of the bimodal hybrid system $\hs$ from Figure \ref{fig:hybrid_quotient_space}, $\hse$, is depicted on the left in Figure~\ref{fig:relaxed_hybrid_quotient_space}. 

Concretely, for each $e\in \Gamma$ we define the \emph{relaxed strip}
$S_e^\epsilon = \set{p + \hat{g}_eq \in \R^n \colon p \in G_e \text{ and } q \in \sp{0, \epsilon}}$
 and then for each $j \in \J$ define the \emph{relaxed domain} $D_j^\epsilon =D_j \cup_{e \in \mathcal{N}_j} S_e^\epsilon.$
Next, for each $e =(j,j') \in \Gamma$ we then define the \emph{relaxed guard set}
$
G_e^\epsilon = \set{x \in S_e^\epsilon \colon g_e^\epsilon(x) = \hat{g}_e^Tx - \p{c_e +\epsilon} = 0  }
$and define the \emph{relaxed reset map} $R_e^\epsilon \colon G_e^\epsilon \to \partial D_{j'}$ by
$R_e^\epsilon(x) = R_e(P_e(x))$,
where $P_e \colon \R^n \to \R^n$ is defined by $P_e(x) = x - \hat{g}_eg_e(x)$. Intuitively, $P_e$ projects points onto the plane containing $G_e$, so that $R_e^\epsilon(G_e^\epsilon) = R_e(G_e)$. We now provide our definition of a \emph{relaxed hybrid system}, which we take from~\cite{burden2015metrization}.
\begin{definition} \label{def:relaxed_hybrid_sys}
Let $\mathcal{H}$ be a hybrid system. We then define the $\epsilon$-relaxation of $\mathcal{H}$ to be the seven-tuple $ 
\mathcal{H}^{\epsilon} =  \p{\mathcal{J}, \Gamma, \mathcal{D}^{\epsilon}, U, \mathcal{F}^{\epsilon}, \mathcal{G}^\epsilon, \mathcal{R}^{\epsilon}}$, where:
\begin{enumerate}
\item $\mathcal{D}^{\epsilon} =  \set{D_j^{\epsilon}}_{j \in \mathcal{J}}$ is the set of relaxed domains ;
\item $\mathcal{F}^\epsilon = \set{f_j^{\epsilon}}_{j \in \mathcal{J}}$ is the set of relaxed vector fields, where $f_j^\epsilon \colon D_j^\epsilon \times U \to \R^n$;
\item $\mathcal{G}^\epsilon = \set{G_e^\epsilon}_{e = \in \Gamma}$ is the set  of relaxed guard set; and,
\item $\mathcal{R}^{\epsilon} = \set{R_e^\epsilon}_{e = \in \Gamma}$ is the set of relaxed reset maps.
\end{enumerate}
\end{definition}

We embed our relaxed hybrid systems in the disjoint union $\coprod_{j \in \J}D_j^\epsilon$ and we adopt the \emph{relaxed hybrid quotient space} introduced in \cite{burden2015metrization}. Let $\hse$ be a relaxed hybrid system and let 
$\widehat{R}^\epsilon \colon \coprod_{e \in \Gamma} G_e^\epsilon \to \coprod_{j \in \J} \partial D_{j}$
be characterized by $\widehat{R}_e^\epsilon(x) = R_e^\epsilon(x)$ for each $x \in G_e^\epsilon$. We then define the relaxed hybrid quotient space of $\hse$ to be
$\Me = \frac{\coprod_{j \in J} D_j^\epsilon }{\Lambda_{\widehat{R}^\epsilon}}$. 

The construction of $\Me$ for our example bimodal hybrid system is shown on the right in Figure \ref{fig:relaxed_hybrid_quotient_space}. Note that the strips $S_{(1,2)}^\epsilon$ and $S_{(2,1)}^\epsilon$ form a single \emph{hybrid strip} in $\Me$ \footnote{This is not technically true, since for each pair of partner edges we only "glued" $G_e^\epsilon$ to $G_{e'}$ and $G_{e'}^\epsilon$ to $G_{e}$. However, it is notationally cumbersome to "glue" the entire width of $S_e^\epsilon$ to $S_{e'}^\epsilon$, and so we choose to abuse notation here.}, and a similar change of coordinates occurs when traversing $e$ in $\Me$, as was described for the same transition in $\M$. The trajectories of our relaxed hybrid systems will again be continuous on $\Me$; such a trajectory $x^\epsilon$ is depicted in Figure~\ref{fig:relaxed_hybrid_quotient_space}, where we also reproduce the trajectory $x$ from Figure~\ref{fig:hybrid_quotient_space}. By representing $x$ on $\Me$, as in~\cite{burden2015metrization}, we will be able to compare the distance between the trajectories of a hybrid system and its $\epsilon$-relaxation using the $\indmetric{\Me}$ metric, which is defined analogously to how $\indmetric{\M}$ was constructed in \eqref{eq:metric}. In particular, given two trajectories $x, x^\epsilon \colon \sp{0,T} \to \Me$, we will use the metric $\rho^\epsilon\p{x, x^\epsilon} = \sup\p{\indmetric{\Me}\p{x(t), x^\epsilon(t)} | t \in \sp{0,T}}$ from~\cite{burden2015metrization}, to bound the distance between different trajectories.

\section{Representing Discrete Jumps with Switched Systems}\label{sec:switched_systems}

We now demonstrate how to describe a discrete transition of a hybrid system using the execution of a switched system, which will allow us to use Filippov solutions to describe the composition of continuous and discrete dynamics along guard sets. We begin by making the change of coordinates that occurs during a discrete transition explicit for a given edge $e =(j,j')$ with partner $e'$.

Definitively, if $\set{v_e^i}_{i = 1}^{n-1}$ is a basis for the subspace parallel to $G_e$, then $\set{v_e^1 ,\dots, v_e^{(n-1)}, \hat{g}_e}$ is a basis for $\R^n$, and when edge $e$ is traversed this basis is transformed element wise to the basis $\set{A_ev_e^1, \dots , A_e v_e^{(n-1)}, - \hat{g}_{e'}}$, where $\set{A_ev_e^i}_{i = 1}^{(n-1)}$ is a basis for the subspace parallel to $G_{e'}$ (indeed this set is linearly independent since we assumed $A_e$ to be full rank), and $-\hat{g}_{e'}$ is orthogonal to this subspace.

\begin{figure}[t]
  \centering
  \includegraphics[width= .6\columnwidth, height = .95\textheight, keepaspectratio ]{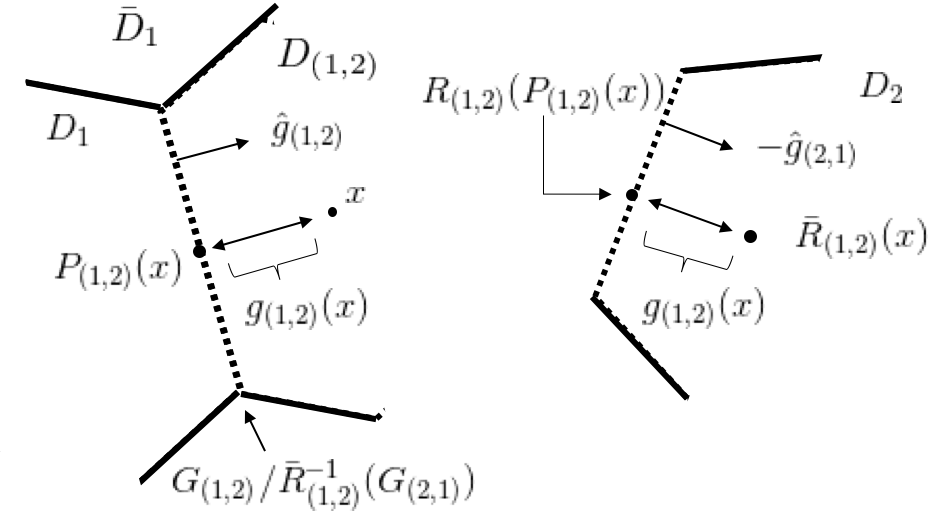}%
  \caption{On the left is depicted $\bar{D}_j = D_j \cup D_{e}$, as well as a point $x \in D_e$. The various components of $\bar{R}_e(x)$ are depicted on $\bar{D}_{j}$ and $D_{j'}$ (on the right).}
  \label{fig:r_bar}
\end{figure}

In order to perform this change of basis automatically during simulation of a discrete transition, we will appropriately translate, rotate and resize  $D_{j'}$, appending it to $D_j$, so we may directly simulate how a trajectory evolves into the interior of $D_{j'}$ after traversing edge $e$.  We denote this transformed version of $D_{j'}$ by $D_{e}$, which is depicted on the left side of Figure~\ref{fig:r_bar}. To accomplish this task define the map $\bar{R}_e \colon \R^n \to \R^n$  by $\bar{R}_e(x) = R_e(P_e(x)) - \hat{g}_{e'}g_e(x)$,
and then define $D_e = \set{x \in \R^n \colon \bar{R}_e(x) \in D_{j'}}$. The  various components of $\bar{R}_e$ are also depicted in Figure~\ref{fig:r_bar}. 

To understand the action of $\bar{R}_e$, recall that $P_e(x)$ projects points onto the plane containing $G_e$, so $R_e(P_e(G_e))= G_{e'}$. Thus, the first term in $\bar{R}_e$ maintains the one-to-one correspondence between $G_e$ and $G_{e'}$ that defines the edge; indeed note that domains $D_j$ and $D_e$ are separated by the surface $G_e \slash \bar{R}_e^{-1}(G_{e'})$. The term $- \hat{g}_{e'}g_e(x)$, on the other hand, aligns the transverse coordinates $\hat{g}_e$ and $- \hat{g}_{e'}$, in the sense that $g_{e'}(\bar{R}_e(x)) = - g_e(x)$, so that as a trajectory $x$ leaves $D_j$ and enters the interior of $D_e$, $\bar{R}_e(x)$ leaves $G_{e'}$ and enters the interior of $D_{j'}$.  Another way to understand $\bar{R}_e$ is to consider the reformulation 
\begin{align}
\bar{R}_e(x) &= R_e(P_e(x)) - \hat{g}_{e'}g_e(x)\\
&= R_e(x - \hat{g}_{e}g_e(x)) - \hat{g}_{e'}g_e(x)\\
&= A_e(I\cdot x -\hat{g}_e(\hat{g}_e^Tx -c_e)) + b_e - \hat{g}_{e'}(\hat{g}_e^Tx - c_e)\\
&= \bar{A}_ex + \bar{b}_e,
\end{align}
where $\bar{A}_e = A_e(I -\hat{g}_e\hat{g}_e^T) - \hat{g}_{e'}\hat{g}_e^T$ and $\bar{b}_e =A_e\hat{g}_ec_e + b_e + \hat{g}_{e'}c_e$.
The matrix $\bar{A}_e$ applies the change of basis that occurs when $e$ is traversed, and is thus invertible. In particular, the matrix $(I -\hat{g}_e\hat{g}_e^T)$ is the natural projection onto the subspace orthogonal to $\hat{g}_e$, so the term $A_e(I -\hat{g}_e\hat{g}_{e}^T)$ applies the change of basis along $G_e$ and $G_{e'}$, while the dyad $-\hat{g}_{e'}\hat{g}_e^T$ rotates vectors in the direction $\hat{g}_e$ to align with $-\hat{g}_{e'}$. These definitions enable the following result.

\begin{lemma}\label{lemma:equiv_flows}
Let $e = (j, j')$. Let $f_e \colon D_e \times U  \to \R^n$ be defined by $f_e(x,u) = \bar{A}_e^{-1} f_{j'}(\bar{R}_e(x),u)$. Then $\forall x \in D_e$ if we take $\deriv{}{t}x = f_e(x,u)$ we have that $\deriv{}{t}\bar{R}_e(x) = f_{j'}(\bar{R}_e(x),u)$.
\end{lemma}

To prove the claim we compute $\deriv{}{t} \bar{R}_e(x)$ = $ \bar{A}_e f_e(x,u)$ = $\bar{A}_e  \bar{A}_e^{-1} f_{j'}(\bar{R}_e(x),u)$=$f_{j'}(\bar{R}_e(x),u)$. Intuitively, $f_e$ evaluates the vector field $f_{j'}$ at the point $\bar{R}_e$, then reverses the change of coordinates that occurs when traversing edge $e$ (by passing the vector field through $\bar{A}_e^{-1}$), effectively transplanting the vector field $f_{j'}$ onto the domain $D_e$.

\begin{figure}[t!]
  \centering
  \includegraphics[width= \columnwidth, height = .95\textheight, keepaspectratio ]{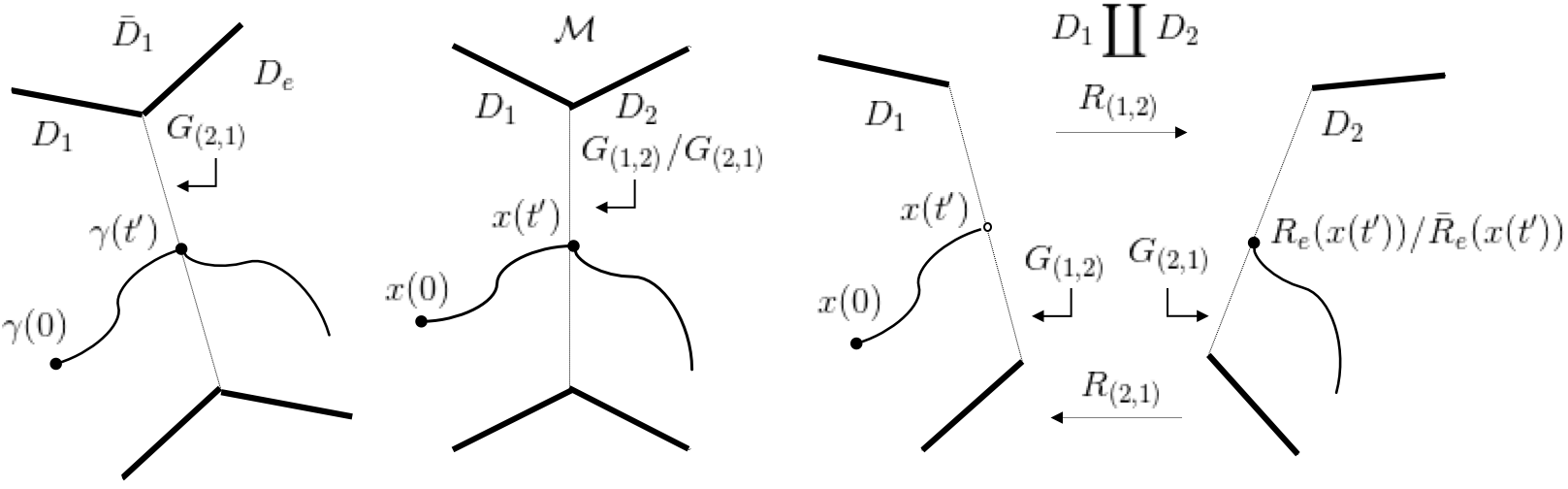}%
  \caption{The curve $\gamma$ is simulated on $\bar{D}_1$, and then interpreted to construct the corresponding transition on $\M$ (center) and $D_1 \coprod D_2$ (right).}
  \label{fig:pure_transition}
\end{figure}

More generally, for a domain $j$ with possibly more than one guard set we define $\bar{D}_j = D_j \cup_{e \in \mathcal{N}_j} D_e$ \footnote{Note that even when two guards $G_e, G_{e'} \subset \partial D_j$ do not intersect, it may be the case that $D_e \cap D_{e'}$ has a non-empty interior. We ignore this technicality, since in practice we will sample neither $f_e$ nor $f_{e'}$ from this region.}, then define the switched system $\bar{f}_j \colon \bar{D}_j \times U \to \R^n$ by
\begin{equation}
\bar{f}_j(x,u)= \begin{cases}
f_j(x,u) & \text{if } x \in D_j\\
f_e(x,u) & \text{if } x \in D_e, \forall e \in \mathcal{N}_j,
\end{cases}
\end{equation}
where $f_e$ is defined as in Lemma~\ref{lemma:equiv_flows}, for each $e \in \mathcal{N}_j$. Using this switched system, we can accurately describe transitions out of mode $j$. For example, suppose the hybrid system is instantiated with initial condition $x(0) \in D_j$, and evolves under the vector field $f_j$ until time $t'$ where $x(t') \in G_e$, for some $e =(j,j') \in \mathcal{N}_j$. The system is then reset to the point $R_e(x(t')) = \bar{R}_e(x(t'))$, and $x$ then evolves under the influence of $f_{j'}$. 
Alternatively, we can simulate the auxiliary curve defined by $\deriv{}{t} \gamma = \bar{f}_j(\gamma,u)$ with initial condition $\gamma(0) = x(0)$, allowing $\gamma$ to flow into the interior of $D_e$. Note that since they share the same differential equation and initial condition, we will have $x(t) = \gamma(t), \forall t \in  [ 0,t')$. At $t'$, we have $\bar{R}_e(\gamma(t')) = x(t') \in D_{j'}$,  and $\forall t \geq t'$ by Lemma~\ref{lemma:equiv_flows} we have $\deriv{}{t}\bar{R}_e(\gamma) = f_{j'}(\bar{R}_e(\gamma), u)$. Thus, $\forall t \geq t'$, we have $\bar{R}_e(\gamma(t)) = x(t)$, since the two curves share the same initial condition and differential equation. Thus, we can construct the trajectory $x$ by simulating $\gamma$, and interpreting $x(t)= \gamma(t) \in D_j$ for $t \in [0,t')$, and interpreting $x(t) = \bar{R}_e(\gamma(t)) \in D_{j'}$ for $t\geq 0$. These curves are depicted in Figure \ref{fig:pure_transition}.

Ultimately, this process empowers us to describe discrete transitions of hybrid systems using the solution concept of Filippov. For edge $e =(j,j')$, carefully inspecting $f_e$, one can see that $\hat{g}_e^Tf_e(x,u) = - \hat{g}_{e'}^Tf_{j'}(\bar{R}_e(x),u)$, so sliding solutions for $\bar{f}_j$ arise along $G_e$ when $f_j$ points into $G_e$ and $f_{j'}$ points into $G_{e'}$, at corresponding points along the hybrid surface $G_e \slash G_{e'}$. At this point we wish to remark that, while many authors (e.g. \cite{johansson1999regularization}) have discussed the possibility of using Filippov solutions to describe sliding Zeno executions for hybrid systems with jumps, to the best of our knowledge, we provide the first explicit means of doing so. Yet, some hybrid transitions which display Zeno phenomena may reduce to a switched system for which the solution concept of Filippov is undefined, since both vector fields are parallel to their respective guard sets. \footnote{For example, the two numerical examples we consider in Section~\ref{sec:examples} fall into this category.}  Our relaxations, however, will resolve this issue.



\section{Relaxed Vector Fields}\label{sec:regularized_vfields}

While we have gained the ability to describe hybrid transitions using the solution concept of Filippov, such trajectories are difficult to approximate numerically, as they require accurately detecting when the guard sets are crossed, and when sliding solution arise and terminate. In order to add slack to our numerical calculations, we extend the method of Teixeira (see e.g. \cite{llibre2007regularization}) to relax our collection of switched systems.

For each edge  $e = (j ,j')$ we define analogs to $\bar{R}_e$ and $D_e$ for the relaxation $\hse$. In particular, we define $\bar{R}_e^\epsilon \colon \R^n \to \R^n$ by $\bar{R}_e^\epsilon(x) = R_e(P_e(x)) + \hat{g}_{e'}g_e^\epsilon(x)$ and then we define 

\begin{figure}
  \centering
  \includegraphics[width= \columnwidth, height = .95\textheight, keepaspectratio ]{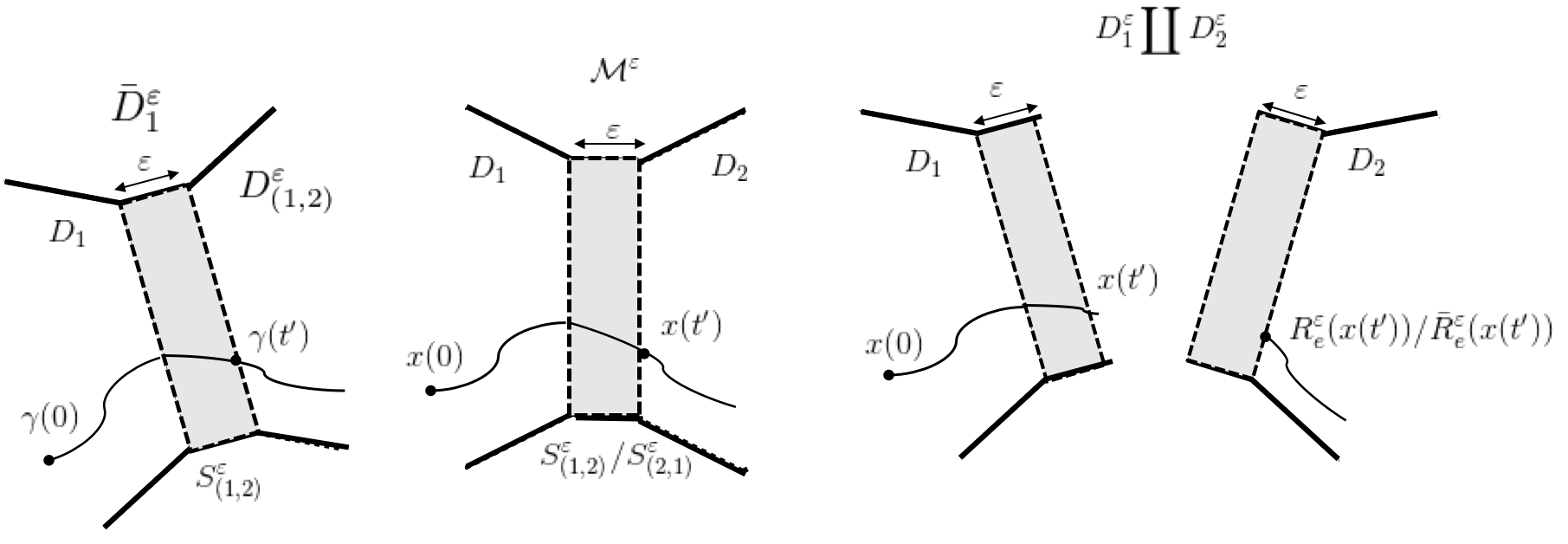}%
  \caption{The curve $\gamma^\epsilon$ is simulated on $\bar{D}_1^\epsilon$, and then interpreted to construct the corresponding transition on $\Me$ (center) and $D_1^\epsilon \coprod D_2^\epsilon$ (right).}
  \label{fig:relaxed_transition}
\end{figure}

\noindent
$D_{e}^\epsilon = \set{x \in \R^n \colon \bar{R}_e^\epsilon(x) \in D_{j'}}$, which is depicted on the left in Figure \ref{fig:relaxed_transition} for our example bimodal hybrid system. Note that $\bar{R}_e^\epsilon(S_e^\epsilon) = S_{e'}^\epsilon$. We may also  refactor $\bar{R}_e^\epsilon(x) = \bar{R}_e(x) + \hat{g}_{e'}\epsilon = \bar{A}_ex + \bar{b}_e + \hat{g}_{e'}\epsilon$. 

We will use the following class of functions from  \cite{llibre2007regularization} to smoothly transition between the dynamics of $f_j$ and and the (projected) dynamics of $f_{j'}$ when crossing $S_e^\epsilon$. We say that $\phi \colon \R \to \sp{0, 1}$ is a \emph{transition function} if $i)$  $\phi(a) = 0$ for $a \leq 0$ and $\phi(a) = 1$ for $a \geq 1$, $ii)$ $\phi'(a) > 0$ for $a \in (0,1)$,  $iii)$ $\phi'$ is Lipschitz continuous, and $iv)$ $\forall a, \phi(1- a)= \phi(a)$ (i.e. $\phi$ is symmetric around 0.5). For the rest of the paper we assume a single transition function has been chosen. \footnote{For example, in our code we employ \begin{equation}\phi(a) =  \begin{cases} 0 & \text{if } a \leq0 \\ \frac{1}{2} + \frac{1}{2}sine(\pi a - \frac{\pi}{2}) & \text{if } 0 < a < 1 \\ 1 & \text{if } 1 \leq a \end{cases}\end{equation}} For the edge $e$ we then define $\phi_e^\epsilon(x) = \phi(\frac{g_e(x)}{\epsilon})$, and now define our relaxation of $f_e$.

\begin{lemma}\label{equiv_flows2}
Let $e =(j ,j')$. Define $f_e \colon D_{j}^\epsilon \cup D_e^\epsilon \times U \to \R^n$ by  $f_e^\epsilon(x,u) = (1-\phi_e^\epsilon(x))f_j(x,u) +  \phi_e^\epsilon(x) \bar{A}_e^{-1}f_{j'}(\bar{R}_e^\epsilon(x),u)$. Then $\forall x \in D_e^\epsilon$ if we take $\deriv{}{t}x = f_e^\epsilon(x,u)$ we have that $\deriv{}{t}\bar{R}_e^\epsilon(x) = f_{j}(\bar{R}_e^\epsilon(x),u)$. \footnote{Again, we compute $\deriv{}{t}\bar{R}_e^\epsilon(x) = \bar{A}_ef_e^\epsilon(x,u) = \bar{A}_e\p{\bar{A}_e}^{-1}f_{j'}(\bar{R}_e^\epsilon(x),u) = f_{j'}(\bar{R}_e^\epsilon(x),u)$.}
\end{lemma}
It was shown in \cite{llibre2008sliding} that for each $\epsilon >0$ the vector field $f_e^\epsilon$ is continuously differentiable. Note that when $g_e(x) \leq 0 $ (and $x \in D_j$),  $\phi_e^\epsilon(x) = 0$ and $f_e^\epsilon(x,u)$ returns $f_j(x,u)$. Similarly, when $g_e(x) \geq \epsilon$ (and $x \in D_{e}^\epsilon$), $f_e^\epsilon$ returns  $\bar{A}_e^{-1}f_{j'}^\epsilon(\bar{R}_e^\epsilon(x),u)$. When $ 0 \leq g_e(x) \leq \epsilon$, $f_j^\epsilon$ produces a convex combination of these vector fields. In the case that $f_j$ points into $G_e$ and $f_{j'}$ points into $G_{e'}$, the trajectories of $f_e^\epsilon$ will remain confined to $S_e^\epsilon$; thus, Zeno executions are approximated by well defined trajectories on our relaxed strips. 

We can use $\bar{R}_e^\epsilon$ and $f_e^\epsilon$ to keep track of how a relaxed trajectory evolves in $D_{j'}$ after a relaxed transition along edge $e = (j,j')$ occurs, employing the same procedure that was developed using $f_e$ and $D_e$ in the previous section. In particular, we simulate the auxiliary curve $\deriv{}{t}\gamma^\epsilon = f_e^\epsilon(\gamma^\epsilon,u)$, allowing this curve to flow through $S_e^\epsilon$ and into $D_e^\epsilon$. We can then use the map $\bar{R}_e^\epsilon$ to keep track of how such a trajectory would have propagated into $D_{j'}$ after crossing $G_e^\epsilon$ and being reset to this domain. This process is depicted in Figure~\ref{fig:relaxed_transition}. Finally, we define $\bar{D}_j^\epsilon =D_j^\epsilon \cup_{e \in \mathcal{N}_j} D_e^\epsilon$ and then define $f_j^\epsilon \colon \bar{D}_j^\epsilon \times U \to \R^n$ by
\begin{equation}
f_j^\epsilon = \begin{cases}
f_j(x,u) & \text{if } x \in D_j \\
f_e^\epsilon(x,u) & \text{if } x \in D_e^\epsilon \bigcup S_e^\epsilon, \ \forall e \in \mathcal{N}_j,
\end{cases}
\end{equation}
which can also be shown to be continuously differentiable, and thus has a Lipschitz continuous gradients, since continuous functions are Lipschitz on compact domains -- a property that will be useful later. For a given relaxed hybrid system $\hse$, we endow the relaxed domain $D_j^\epsilon$ with the vector field $f_j^\epsilon |_{D_j^\epsilon \times U}$, for the purposes of Definition~\ref{def:relaxed_hybrid_system}; however, our theoretical analysis and discrete approximations will rely on simulating $f_j^\epsilon$ past the relaxed strips $\set{S_e^\epsilon}_{e \in \mathcal{N}_j}$ and into the projected domains $\set{D_e^\epsilon}_{e \in \Gamma}$. 

\begin{lemma}\label{lemma:equiv3}
Let $e$ and $e'$ be partner edges. Then $\forall x \in S_e^\epsilon$, if we take $\deriv{}{t}x = f_e^\epsilon(x,u)$, then $\deriv{}{t}\bar{R}_e^\epsilon(x) = f_{e'}^\epsilon(\bar{R}_e^\epsilon(x),u)$. 
\end{lemma}

In other words, the vector fields $f_e^\epsilon$ and $f_{e'}^\epsilon$ produce \emph{equivalent flows} over $\Me$, and thus we can represent a relaxed transition along $e =(j,j')$ on either $S_e^\epsilon$ or $S_{e'}^\epsilon$. Moreover, this implies that, if a relaxed trajectory repeatedly flows back and forth across $S_e^\epsilon \slash S_{e'}^\epsilon$, we can simulate this behavior on either $\bar{D}_j^\epsilon$ or $\bar{D}_{j'}^\epsilon$, and don't need to switch between the two vector fields each time a transition occurs -- this fact will greatly simplify out analysis later. The proof of Lemma \ref{lemma:equiv3} is largely algebraic, and uses the fact that $\phi$ is symmetric about 0.5. 

We conclude this section by studying how the trajectories of $f_j^\epsilon$ converge to those of $\bar{f}_j$ as we take $\epsilon \to 0$. For the following two theorems, assume we have fixed an input $u \in \BVU$, and then let $x^\epsilon \colon \sp{0,T} \to \bar{D}_j^\epsilon$ be the corresponding solution generated by $f_j^\epsilon$ with initial condition $x_0^\epsilon \in D_j$, and let $\bar{x} \colon \sp{0,T} \to \bar{D}_j$ be the trajectory generated by $\bar{f}_j$ with initial condition $x_0 \in D_j$ and the same input. We leave the proofs to the appendix. A version of the following result for autonomous vector fields may be found in \cite{fiore2016contraction}.

\begin{theorem}\label{thm:regular_convergence}
Assume that for each $e =(j,j') \in \mathcal{N}_j$ and each $(x,u) \in G_e \times U$ either $\hat{g}_e^Tf_{j}(x,u) >0$ or $\hat{g}_e^Tf_{e}(x,u) < 0$. Finally assume $\norm{x^\epsilon_0 - x_0} \leq k\epsilon$, for some $k>0$. \footnote{During each discrete transition our relaxations will incur an error of order $\epsilon$. By adding error to our initial conditions here, we will be able to call this result inductively to prove convergence when trajectories undergo multiple transitions.} Then $\exists \epsilon_0 > 0$ and $C >0$ such that for each $\epsilon_0 >\epsilon > 0$ we may bound $\norm{\bar{x} - x^\epsilon}_{\infty} \leq C \epsilon$.
\end{theorem}

The hypothesis of Theorem~\ref{thm:regular_convergence} guarantee that Filippov solutions are unique and well defined for $\bar{f}_j$ along the guard sets of $D_j$, since the escaping region is empty. Thus, our relaxed vector fields converge to Filippov solutions, when applicable. We next examine how our relaxations behave when Filippov solutions are ill-defined.

\begin{theorem}\label{thm:irregular_convergence}
Suppose there exists a Lipschitz continuous function $l \colon \R \to \R^n$ such that $ l(\epsilon) = x_0^\epsilon $ and $l(0) = x_0$.\footnote{Again we add slack to our initial condition so this result may be called inductively.} Then there exists a uniformly continuous $x^0 \colon \sp{0,T} \to \bar{D}_j$ such that $x^0(0) = x_0$  and $\norm{x^0 - x^\epsilon}_{\infty} \to 0$ as as $\epsilon \to 0$.
\end{theorem}
Theorem \ref{thm:irregular_convergence} implies that our relaxations converge uniformly to a unique, well defined limit, even when the solution concept of Filippov breaks down.

\section{Executions}\label{sec:executions}
Having demonstrated our relaxation approach to describe single discrete transitions of a hybrid system, we modify the algorithmic construction presented in \cite{burden2015metrization} to define the trajectories of our relaxations through multiple transitions.
\begin{definition}\label{def:relaxed_execution} 
An execution for a relaxed hybrid dynamical system $\hse$, given data $x_0 \in D_j$ and $u \in \BVU$, denoted $x^\epsilon \colon \sp{0,T} \to \M^\epsilon$ is constructed via the following algorithm.
\begin{enumerate}
\item Set $x^\epsilon(0)= x_0$ and $t=0$, and let $j \in \J$.
\item Simulate the differential equation $\dot{\gamma^\epsilon}(s) = f_j^\epsilon(\gamma^\epsilon(s),u(s))$ forward in time with initial condition $\gamma^\epsilon(t) = x^\epsilon(t)$ until time $t' = \min\set{T, \inf\set{s \colon \gamma^\epsilon(s^+) \notin  D_j^\epsilon}}$.
\item If $t' = T$ or $\nexists e \in \mathcal{N}_j$ such that $\gamma^\epsilon(t') \in G_e^\epsilon$, let $x^\epsilon(s) = \gamma^\epsilon(s), \ \forall s \in \sp{t,t'}$. Then terminate the execution. 
\item Else let $e = (j,j')$ be such that $\gamma^\epsilon(t') \in G_e^\epsilon$. For each $s \in [t ,t')$ set $x^\epsilon(s) = \gamma^\epsilon(s)$. Set $x^\epsilon(t') = R_e^\epsilon(\gamma^\epsilon(t'))$, set $t =t'$ and set $j = j'$. Go to step 2.
\end{enumerate}
\end{definition}
First, we note that the only time the execution terminates is in line 3, when either the simulation horizon $T$ has been reached or when the trajectory leaves a relaxed domain at a point that does not belong to a relaxed guard set.
Second, we note that the trajectories generated in Definition \ref{def:relaxed_execution} agree with typical definitions for the execution of a hybrid system; that is, a differential equation is simulated until a guard is reached, then the state is reset and resumes simulation. However, these trajectories are continuous over $\Me$ \cite{burden2015metrization}, and are even Lipschitz continuous with respect to their arguments.

\begin{proposition}
Construct $x_1^\epsilon,x_2^\epsilon \colon \sp{0,T} \to \Me$ as in Definition~\ref{def:relaxed_execution} using the arguments $(x_0^1, u^1), (x_0^2, u^2) \in \M^\epsilon \times \BVU$, respectively. Then there exists $L>0$ such that $\rho^\epsilon(x_1, x_2) \leq L\p{\indmetric{\Me}(x_0^1, x_0^2) + \norm{u^1 - u^2}_2}$. 
\end{proposition}
Note that this fundamental systems theoretic property is missing from previous relaxation approaches as \cite{burden2015metrization} and \cite{johansson1999regularization}. Due to space constraints, we do not formally compute variations over our relaxations. However, the result follows from two observations. Firstly, as demonstrated by Lemmas~\ref{equiv_flows2} and \ref{lemma:equiv3}, each portion of a relaxed execution constructed via \ref{def:relaxed_execution} has a one-to-one, affine (and therefore Lipschitz) correspondence to the trajectories generated by a vector field that has Lipschitz continuous gradients. Secondly, by Theorem 5.6.7 of~\cite{polak2012optimization}, the flows generated by each of these vector fields are Lipschitz continuous with respect to their arguments. In Section~\ref{sec:examples} we demonstrate how to compute variations through a relaxed transition in the numerical setting. Although we must construct their trajectories in an algorithmic manner, our class of relaxed hybrid systems may largely be viewed simply as classical dynamical systems -- that is, systems whose trajectories are continuous and have variations which are Lipschitz continuous. Furthermore, the convergence results of Theorems~\ref{thm:regular_convergence} and \ref{thm:irregular_convergence} hold when multiple transitions occur. Note that the following construction is similar to the definition of an execution of a hybrid system from \cite{burden2015metrization}, but unlike this work we are able to describe sliding solutions along our guard sets.
 
\begin{theorem}\label{thm:convergence_1}
Assume that for each $e =(j,j') \in \Gamma$ and each $\p{x,u} \in G_e \times U$ either $\hat{g}_e^Tf_j\p{x,u} > 0$ or $\hat{g}_e^Tf_e(x,u) < 0$. For each $x_0 \in D_j$ and $u \in \BVU$ let $x^\epsilon$ be constructed via Definition \ref{def:relaxed_execution}. Then $\exists \epsilon_0 > 0$ and $C>0$ such that $\forall \epsilon_0 > \epsilon > 0$, $\rho^\epsilon\p{x, x^\epsilon} \leq C \epsilon$, where $x \colon \sp{0,T} \to \M$ is generated by the following algorithm. \begin{enumerate}
\item Set $x(0)= x_0$ and $t=0$, and let $j \in \J$.
\item Simulate the differential equation $\dot{\gamma}(s) = \bar{f}_j(\gamma(s),u(s))$ forward in time (using the solution concept of Filippov) with initial condition $\gamma(t) = x(t)$ until time \\ $t' = \min\set{T, \inf\set{s \colon \gamma(s^+) \notin D_j}}$.
\item If $t' = T$ or $\nexists e \in \mathcal{N}_j$ such that $\gamma(t') \in G_e$, let $x(s) = \gamma(s), \ \forall s \in \sp{t,t'}$. Then terminate the execution.
\item Else let $e = (j,j')$ be such that $\gamma(t') \in G_e$. For each $s \in [t ,t')$ set $x(s) = \gamma(s)$. Set $x(t') = R_e^\epsilon(\gamma(t'))$, set $t =t'$ and set $j = j'$. Go to step 2.\end{enumerate}
\end{theorem}
We again leave the proof to Appendix~\ref{sec:proofs}.
\begin{theorem}\label{thm:convergence_2}
Fix $x \in D_j$, $u \in \BVU$ and let $x^\epsilon$ be constructed by the algorithm in Definition \ref{def:relaxed_execution}. Then
there exists a uniformly continuous $x \colon \sp{0,T} \to \M$ such that $\rho^\epsilon\p{x, x^\epsilon} \to 0$ as $\epsilon \to 0$, where $x(0) =x_0 \in D_j$. 
\end{theorem}
We omit the proof in the interest of brevity since the proof is analogous to that of Theorem~\ref{thm:convergence_1}, except Theorem~\ref{thm:irregular_convergence} is called inductively in place of \ref{thm:regular_convergence}. We employ this limit to define the execution of our relaxed hybrid systems. 
\begin{definition}\label{def:pure_execution}         
Let $\hs$ be a hybrid system. Then given data $x_0 \in \M$ and $u \in \BVU$ we define the corresponding trajectory of $\hs$ to be $x \colon \sp{0,T} \to \M$, where $x = \lim_{\epsilon \to 0} x^\epsilon$, and for each $\epsilon > 0$ we construct $x^\epsilon$ using the algorithm in definition \ref{def:relaxed_execution}. 
\end{definition}

Taken together, Theorems~\ref{thm:convergence_1} and~\ref{thm:convergence_2} imply that executions of our hybrid systems, as in Definition \ref{def:pure_execution}, are unique and we defined, even when traditional solution concepts for hybrid systems would have produced Zeno executions. Note that his property is fundamental, yet missing from current methods such as \cite{johansson1999regularization}, \cite{burden2015metrization}, and \cite{1656623}. While further work is needed to more carefully characterize this limit in cases where Filippov solutions are ill-defined, in Section~\ref{sec:examples} we provide numerical evidence that in such cases our relaxations converge to solutions which make physical sense. 

We now introduce the provably convergent numerical integration scheme that we use to approximate the trajectories of our relaxed hybrid systems. Again, our discretization scheme is largely similar to the one proposed in~\cite{burden2015metrization}. We begin with the following definition of a numerical integrator. 

\begin{definition}\cite{burden2015metrization}
Given a relaxed hybrid system $\hse$, we say $\mathcal{A} \colon \R^n \times U \times \J \times \R \to \R^n$ is a numerical integrator of order $\omega$, if for each $j \in \J$ and $h = T/N$ (where $N \in \N$), and each $x_0 \in D_j$ and $u \in \BVU$ we have 
\[ \sup\p{\norm{x(kh)- z^{\epsilon,h}(kh)} \colon k \in \set{0 , 1, \dots , N}} = O(h^\omega), \]
where $x(0) = x_0$ and $\deriv{}{t}x = f_j^\epsilon(x,u)$, and $z^{\epsilon, h}(0)$ and \\ $z^{\epsilon,h}((k+1)h) = \mathcal{A}(z(kh),u(kh), j ,h)$. 
\end{definition}
As was noted in~\cite{burden2015metrization}, this definition of a numerical integrator is compatible with a large class of discretization schemes, including Euler and the Runge-Kutta family. 

\begin{figure}
  \centering
  \includegraphics[width= \columnwidth, height = .95\textheight, keepaspectratio ]{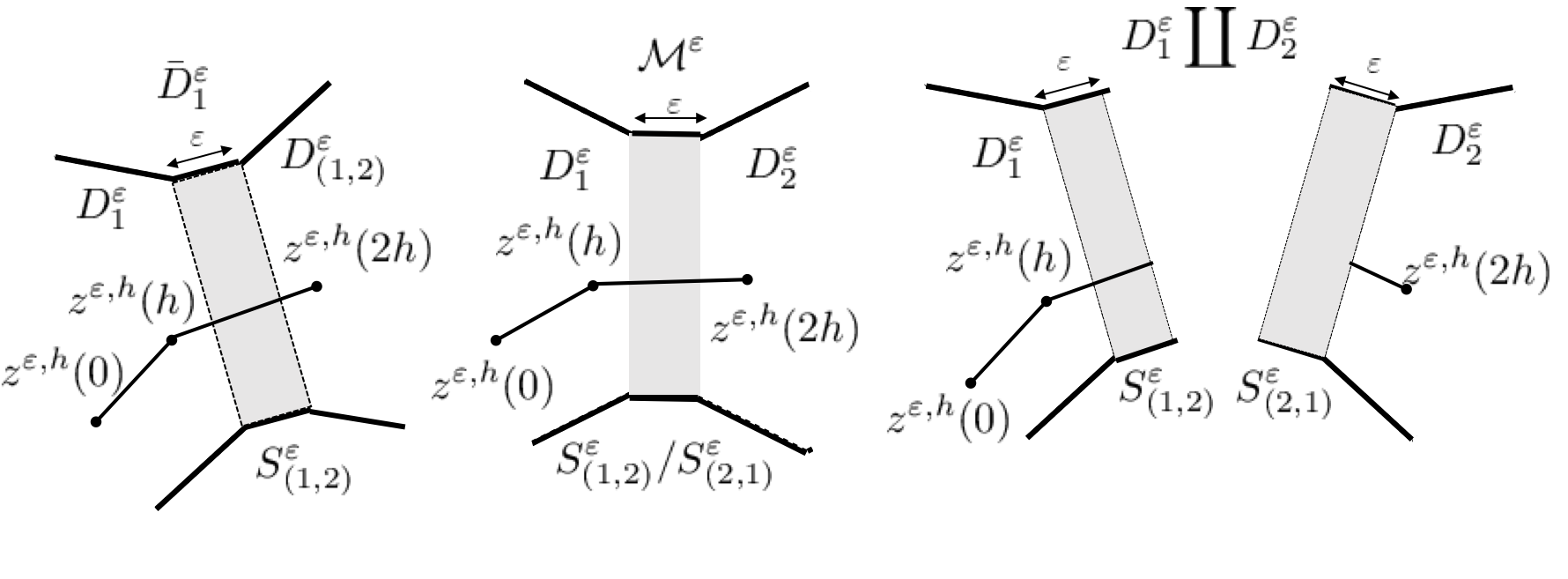}%
  \caption{A numerical approximation of a discrete transition is constructed on $\bar{D}_1^\epsilon$ (left) and the interpreted on $\Me$ (center) and $D_1^\epsilon \coprod D_2^\epsilon$ (right).}
  \label{fig:d_approx}
\end{figure}

\begin{definition}\label{def:d_approx}
Given a relaxed hybrid system $\hse$, initial condition $x_0 \in D_j$, input $u \in \BVU$, step size $h = \frac{T}{N}$ (where $N \in \N$), we construct the discrete approximation $z^{\epsilon,h} \colon \sp{0,t} \to \Me$ according to the following algorithm. 
\begin{enumerate}
\item Let $z^{\epsilon,h}(0) = x_0$, $t =0$, $k =0$ and $j \in \J$.
\item If $k = N$, terminate the execution. Otherwise, let $\gamma^{k+1} = \mathcal{A}\p{z^{\epsilon,k}(kh), u(kh), j, h}$.
\item For each $t \in [kh ,(k+1)h)$ set  \\
$z^{\epsilon, h}(t) = \frac{(k+1)h - t}{h}\gamma^{k+1} + \frac{t -kh}{h}z^{\epsilon, h}(kh)$. 
\item If $\gamma^{k+1} \notin \bar{D}_j^\epsilon$, then let $\bar{t} = \inf\set{t \colon z^{\epsilon, h}(t) \in \bar{D}_j^\epsilon}$ and return $z^{\epsilon, h}|_{\sp{0 ,\bar{t}}}$. Terminate the execution.
\item If $ \exists e =(j,j') \in \mathcal{N}_j$ such that $\gamma^{k+1} \in D_e^\epsilon$, set \\ $ z^{\epsilon, h}((k+1)h)= \bar{R}_e(\gamma^{k+1})$, set k = k+1, and set $j = j'$. Go to step 2.
\item Otherwise, set $z^{\epsilon,h}((k+1)h)$ and k = k+1. Go to step 2. 
\end{enumerate}
\end{definition}  
Our definition of a numerical approximation for relaxed hybrid systems differs from~\cite{burden2015metrization} in one crucial way. The discretization scheme proposed in~\cite{burden2015metrization} requires that a time step be placed in a relaxed strip when simulating a discrete transition. This requires many sample steps be taken until one is placed correctly. On the other hand, our numerical approximation can step \emph{over} the relaxed strip $S_e^\epsilon$ when approximating a discrete transition along edge $e$, since we can use the Lipschitz vector field $f_e^\epsilon$ and map $\bar{R}_e^\epsilon$ to simultaneously model how the trajectory evolves on either side of the transition, and thus there is no need to modulate the step size during numerical approximation of a discrete transition, as is depicted in Figure~\ref{fig:d_approx}.

However, the proof of convergence for our algorithm follows directly from an argument similar to that of Theorem 27 of~\cite{burden2015metrization}. In particular, since, for each $\epsilon > 0$ the vector fields we integrate over are Lipschitz continuous, we incur a numerical error of order $O(h^\omega)$ on each mode, and using an argument similar to that of Theorem 27 of~\cite{burden2015metrization} \footnote{Alternatively, one can make an argument similar to the proof of Theorem~\ref{thm:convergence_1}.}, one can show for each $x_0 \in D_j$, $u \in \BVU$, and step size $h$ small enough that $\rho^\epsilon\p{x^\epsilon, z^{\epsilon, h}} \leq C h^\omega$ for some $C>0$, were $x^\epsilon$ and $z^{\epsilon,h}$ are constructed via Definitions~\ref{def:relaxed_execution} and \ref{def:d_approx}, respectively. Using these conditions, if we construct the hybrid executions $x$ using Definition~\ref{def:pure_execution}, applying the triangle inequality on $\rho^\epsilon$, one can then show that $ \lim_{\epsilon \to 0} \lim_{h \to 0} \rho^\epsilon\p{x, z^{\epsilon,h}} = 0$, as was demonstrated in Corollary 28 of \cite{burden2015metrization}. Moreover, the rate of convergence in $h$ is of order $\omega$. When the hypothesis of Theorem~\ref{thm:convergence_1} are satisfied, the rate of convergence is linear in $\epsilon$, but unknown when Filippov solutions are ill-defined on the guard sets. Before proceeding to our examples, we note that the relaxation scheme we developed in this paper allowed us to construct a provably convergent numerical algorithm capable of simulating all of the trajectories of our hybrid systems, even those that continue past Zeno, an improvement over existing methods such as \cite{burden2015metrization} and \cite{Esposito2001}.

\section{Numerical Examples}\label{sec:examples}
We present several numerical examples, demonstrating the utility of the techniques developed in this paper. We first use an example which is often used to model limbs in the dynamic walking literature. This example is inspired by \cite{burden2011numerical} and \cite{or2009formal}.

\begin{figure*}[h!]
\centering
\begin{tabular}{ c c c}
  \subcaptionbox{}[.33\linewidth]{%
    \includegraphics[width=.35\linewidth]{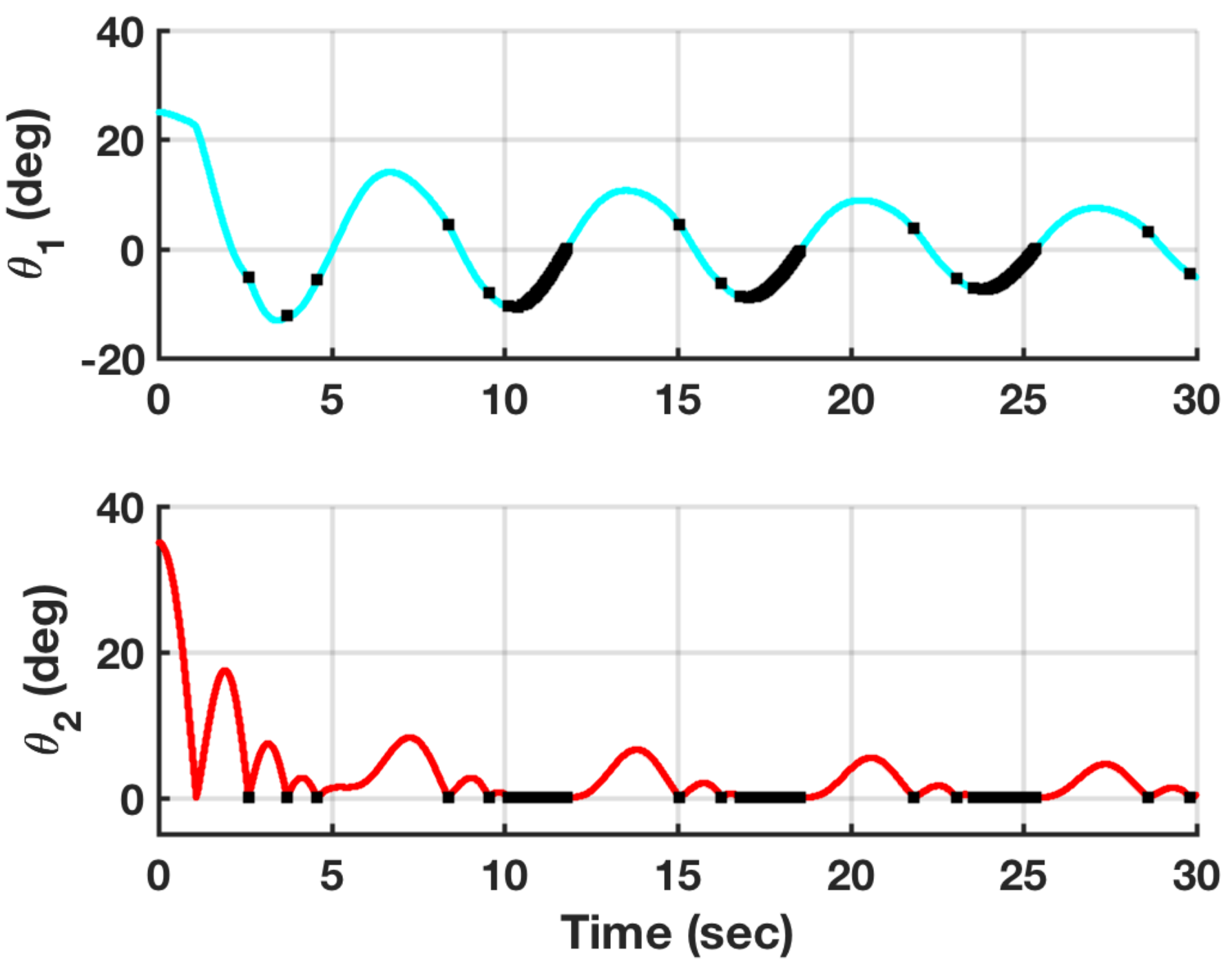}}
    &
  \subcaptionbox{}[.33\linewidth]{%
    \includegraphics[width=.35\linewidth]{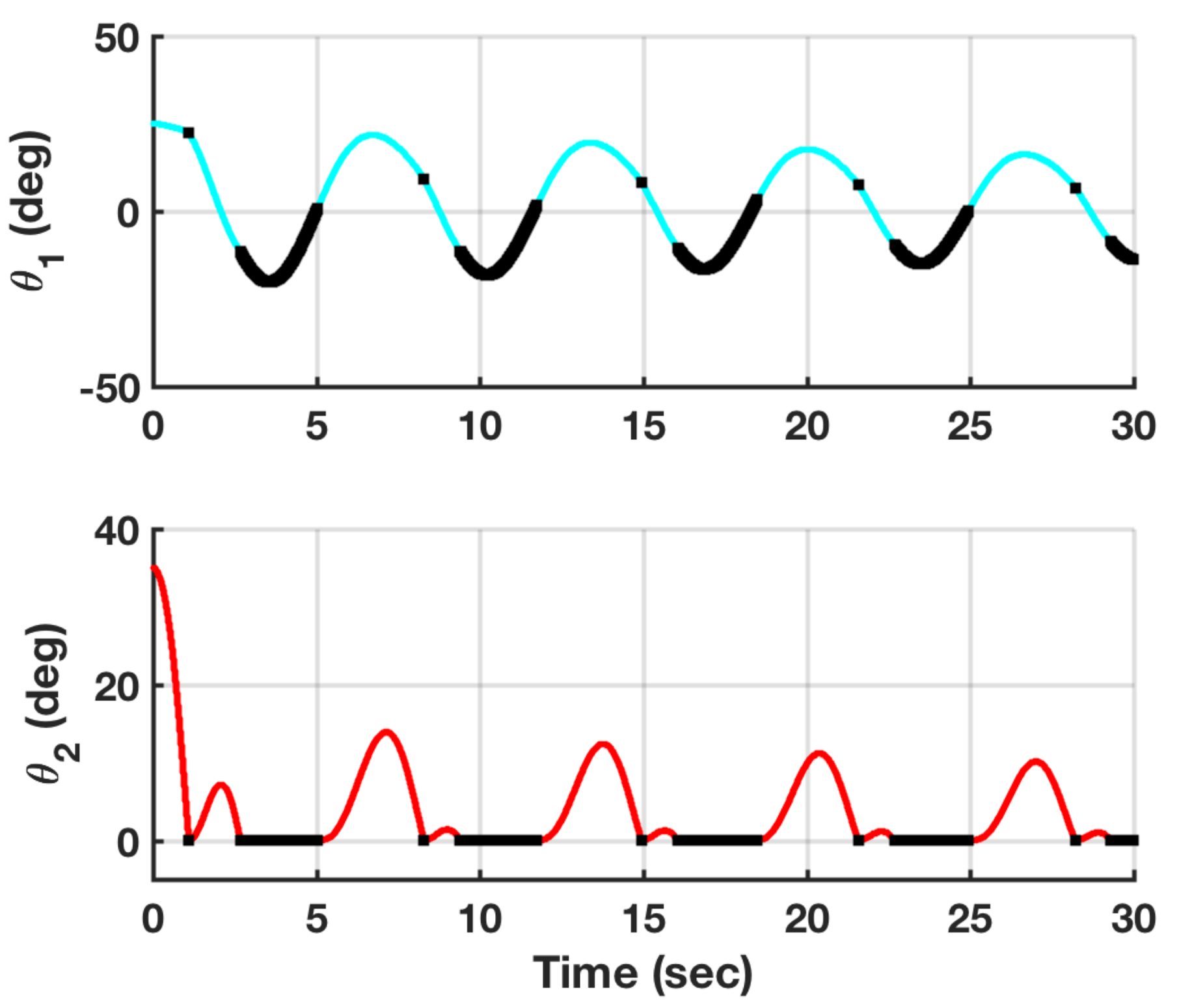}}
    &
  \subcaptionbox{}[.32\linewidth]{%
    \includegraphics[width=.25\linewidth]{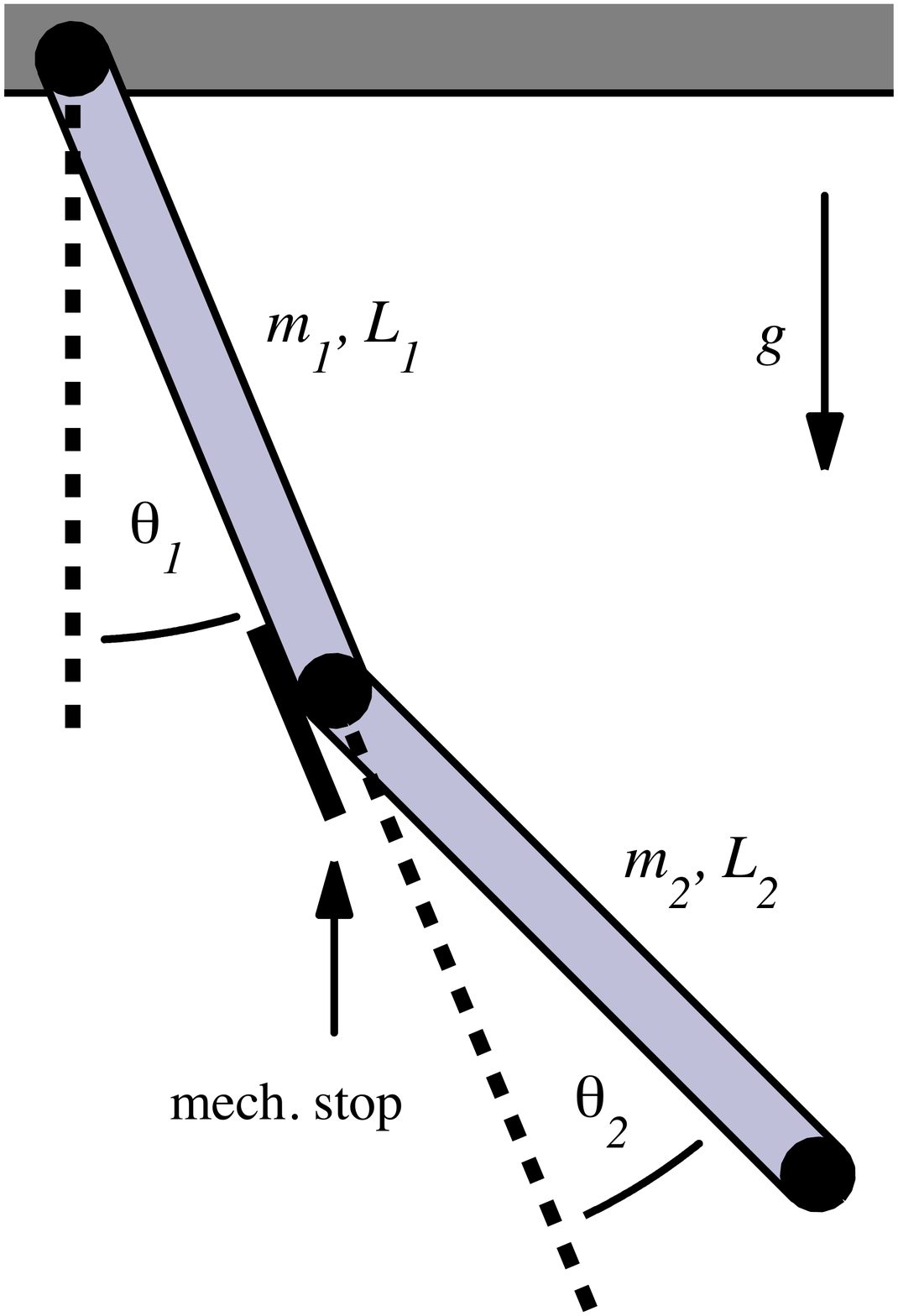}}\\
    \\
    \subcaptionbox{}[.33\linewidth]{%
    \includegraphics[width=.35\linewidth]{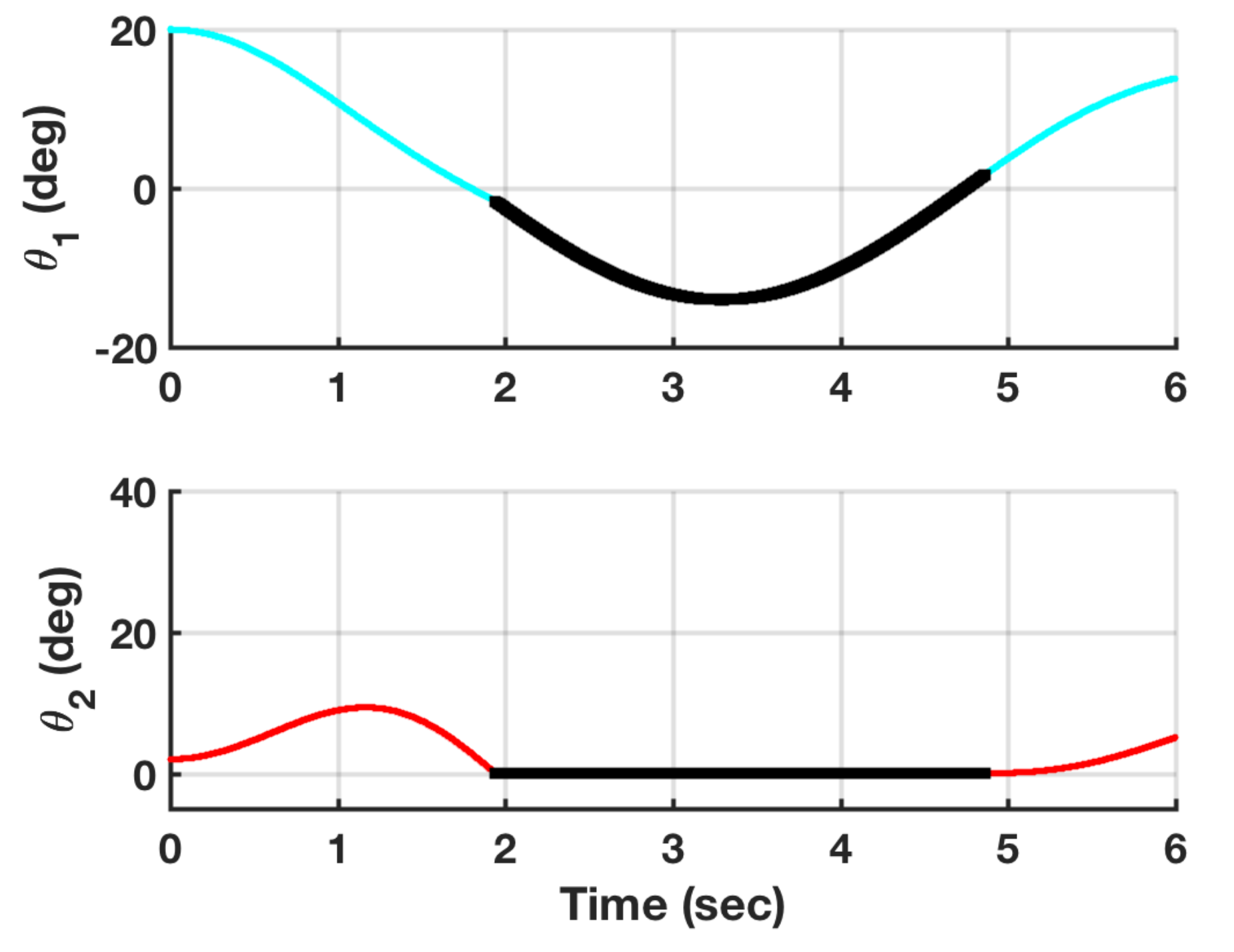}}
    &
    \subcaptionbox{}[.33\linewidth]{%
    \includegraphics[width=.35\linewidth]{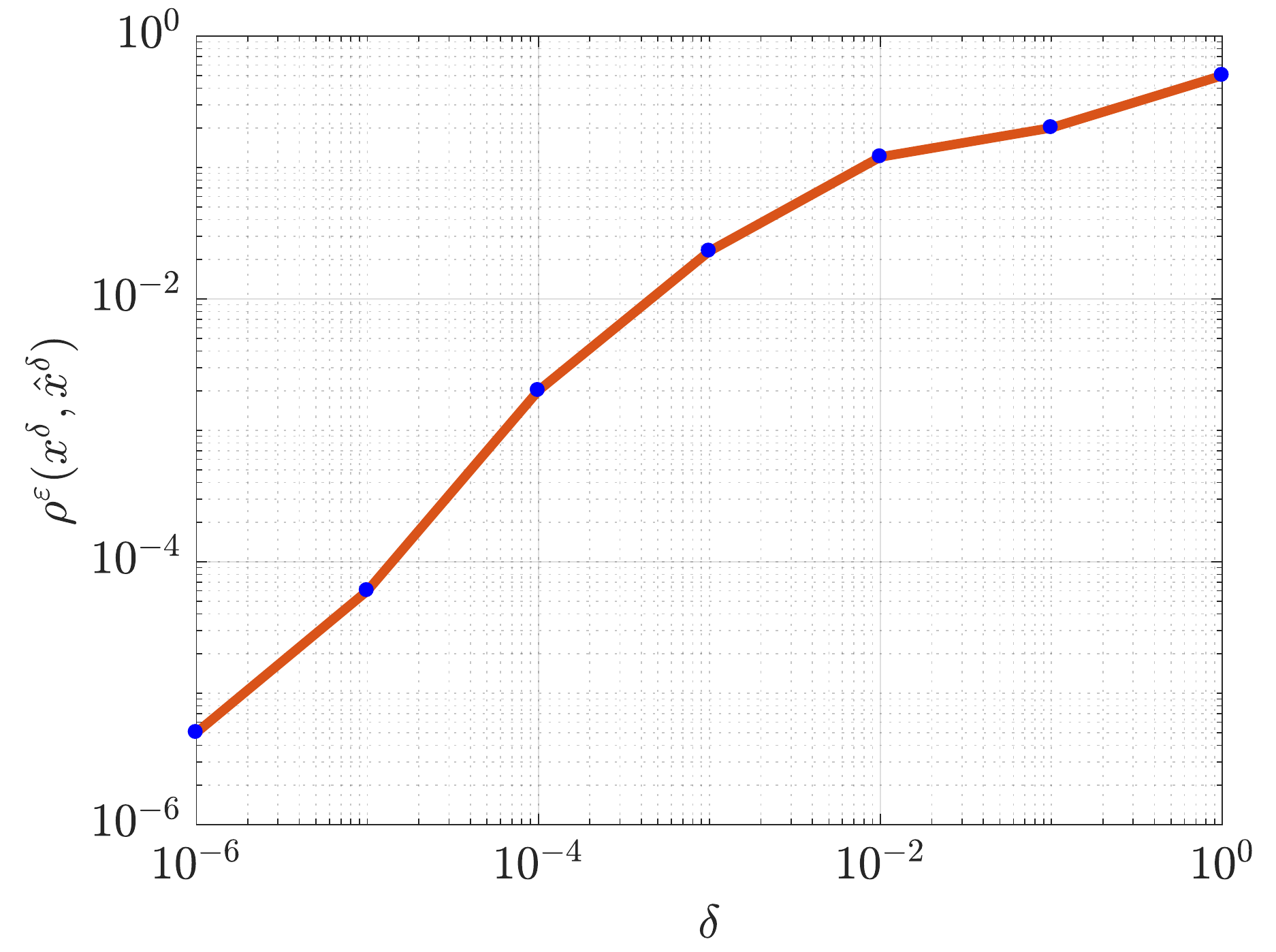}}
    &
    \subcaptionbox{}[.33\linewidth]{%
    \includegraphics[width=.35\linewidth]{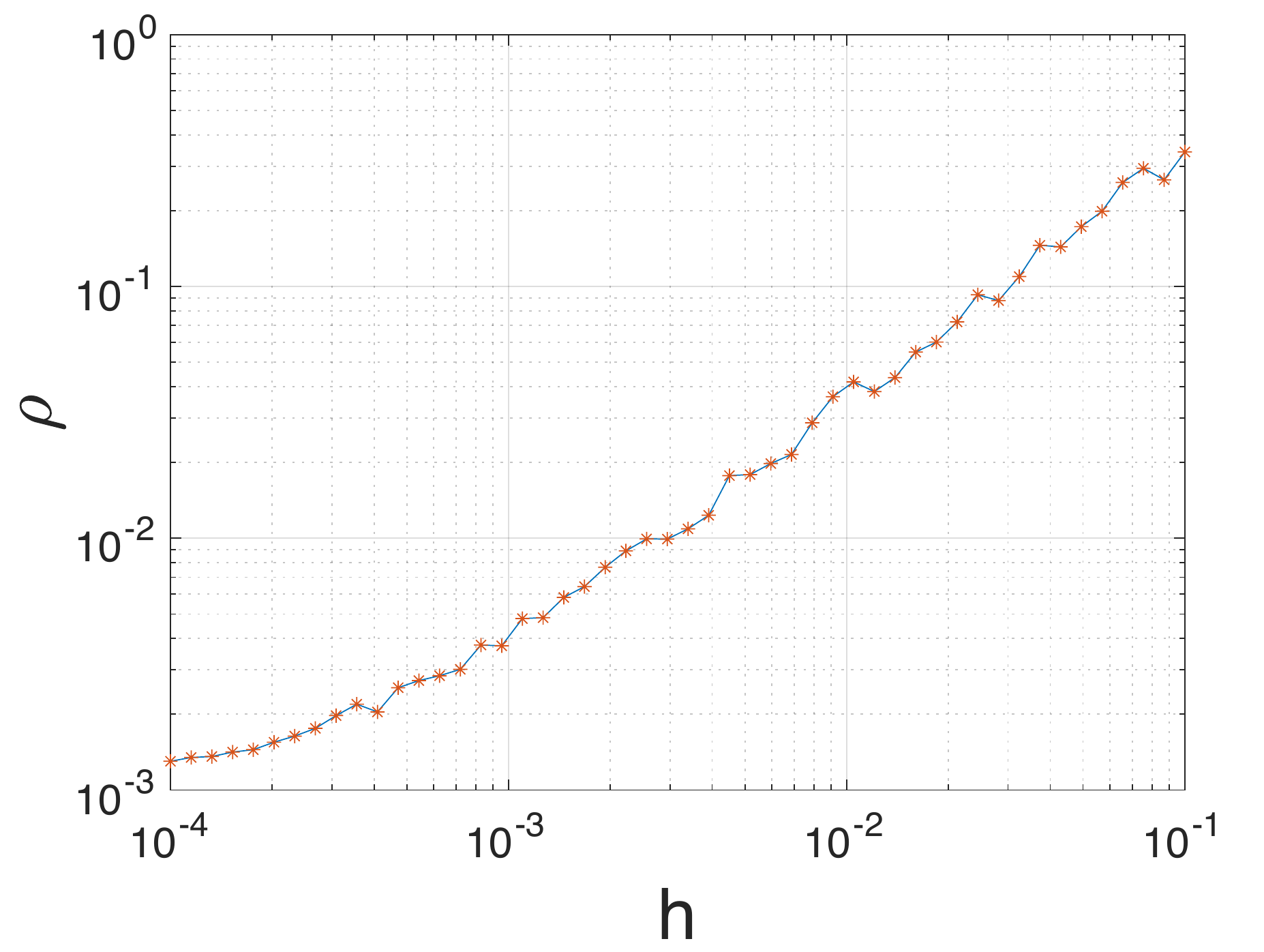}}
    \end{tabular}
    
      \caption{Double Pendulum with a mechanical stop: a) simulation with $c=0.5$ b) simulation with $c = 0$ c) schematic of the system d) reference trajectory for sensitivity analysis e) numerical error of sensitivity analysis. Bouncing Ball: f) simulation error under $\rho$ metric.}\label{fig:dp}
\end{figure*}

\underline{\emph{Example 1: }} (Double Pendulum) Consider the double pendulum with a mechanical stop which is depicted in Figure~\ref{fig:dp}. The system has two angular degrees of freedom $q = (\theta_1, \theta_2)$ whose dynamics are Lagrangian. When the second link impacts the mechanical stop (i.e when $\theta_2 = 0$ and $\dot{\theta}_2 \leq 0$), the angular velocities of are reset according to $(\dot{\theta}_1, \dot{\theta}_2) \to (\dot{\theta}_1+ k(1+c)\dot{\theta}_2, -c\dot{\theta}_2 )$, where $c \in \sp{0,1}$ is the \emph{coefficient of restitution}, and $k > 0$. The interested reader may find the explicit representations of these dynamics in \cite{or2009formal}, where it was demonstrated this system may be faithfully modeled by a unimodal hybrid system with a single edge. When $\theta_2 =\dot{\theta}_2 = 0$, the arm may be locked in place until the imaginary force $\lambda(q, \dot{q})$ becomes nonpositive, at which point the second arm begins to swing freely again. It was shown in \cite{or2009formal} that this hidden locked mode corresponds to a Zeno execution. However, using our relaxation procedure, we can model the dynamics of this hidden mode using well defined solutions on the relaxed strip for this hybrid system. In Figure~\ref{fig:dp}, we simulate trajectories for this system for both $c = 0.5$ and $c=0$, with physical parameters $m_1 = m_2 = L_1 = L_2 = g = 1$, Euler step size $h = 10^{-6}$, $\epsilon = 10^{-6}$ and initial condition $\p{\theta_1, \dot{\theta}_1, \theta_2, \dot{\theta}_2} = \p{25^{\circ}, 0, 35^{\circ}, 0}$ (using the extensions to our framework outlined in the optional appendix). In both simulations, time steps that lie in the relaxed strip are bold and colored black. Note, in both cases the double pendulum settles into a (decaying) periodic orbit, wherein the second arm is periodically locked into place (and the simulation remains confined to the strip) until the imaginary force dissipates and the second arm swings freely. 

For our second experiment involving the double pendulum, we demonstrate how our relaxation framework may be used to conduct sensitivity analysis around a nominal trajectory, as the trajectory progresses through a hybrid transition. In particular, we fix $c =0$ and $\epsilon = 10^{-3}$, and once again fix $m_1 = m_2 = L_1 = L_2 = g = 1$, and choose the initial condition of $\p{\theta_1, \dot{\theta}_1, \theta_2, \dot{\theta}_2}= \p{20^{\circ}, 0, 2^{\circ}, 0} = x_0$ for our nominal trajectory, which is depicted in Figure~\ref{fig:dp} d). We choose again choose an Euler step-size of $h = 10^{-6}$ for each of the simulations of this experiment. Time instances that lie in the relaxed strip are again blackened. Note, this nominal trajectory only undergoes one transition, thus we can simulate the entirety of this trajectory on a single extended domain, without needing to ever reset the trajectory. We denote this domain $D$ and its vector field by $f$. Moreover, since this vector field has gradients that are Lipschitz continuous we can numerically approximate variations over this vector field, using the techniques of, e.g., Chapter 5.6.3 of \cite{polak2012optimization}.   For a given $\delta > 0$ we let $x^\delta$ be the trajectory corresponding to the perturbed initial condition $x_0^\delta= x_0 + \delta x_0$, where $\delta x_0 = \delta(0,0,1^{\circ}, 0)$. Next, applying Theorem 5.6.13 of \cite{polak2012optimization}, and linearizing about the nominal trajectory $x^0$, for each $\delta$ we approximate $x^\delta$ with $\hat{x}^\delta = x^0 + Dx^\delta$ where $Dx^\delta$ is the solution to the linearized difference equation
\begin{equation}
Dx^\delta((k+1)h) = \deriv{}{x}f(x^0(kh))Dx^{\delta}\p{kh}
\end{equation}
with initial condition $Dx^\delta(0) = \delta x_0$. For various values of $\delta$ we simulate both $x^\delta$ and $\hat{x}^\delta$, and in Figure~\ref{fig:dp} e) we appropriately interpret $x^\delta$ and $\hat{x}^\delta$ on $\Me$, and plot the difference $\rho^\epsilon\p{x^\delta, \hat{x}^\delta}$. As this figure clearly demonstrates, using this technique we are able to accurately compute variations through a relaxed transition, even for a trajectory that is simulating past Zeno, as we take $\delta$ to be sufficiently small. 

\underline{\emph{Example 2: }} (Bouncing Ball) For our second example, we simulate the famous bouncing ball. This system consists of a ball repeatedly bouncing on the ground, losing a fraction of its energy during each impact. The ball bounces vertically and has two continuous states, its height $x_1$ and its vertical velocity $x_2$. These two states evolve according to $\deriv{}{t}x_1 = x_2$ and $\deriv{}{t} x_2 = -g$, where $g$ is the gravitational constant. When an impact occurs, the velocity is reset according to $x_2 \to -c x_2$, where $c \in \sp{0,1}$ is again the coefficient of restitution. It can be shown \cite{johansson1999regularization} that the ball undergoes an infinite number of bounces by the finite time $t_{\infty} = \frac{x_2(0)}{g} + \frac{(1/c +1)\sqrt{x_2(0)^2 + 2gx_1(0)}}{g(1/c -1)}$, at which time it comes to rest (i.e. $x_1(t) = x_2(t) =0, \forall  t \geq t_{\infty}$). Thus, a faithful hybrid representation of the system is necessarily Zeno. We simulate this example to benchmark the performance of our relaxations, since we know analytically when and where Zeno occurs. We simulate the bouncing ball with initial condition $(x_1(0), x_2(0)) = (1,0)$, $g =1$ and $c =0.5$ for various Euler steps sizes $h$, and for each simulation fix $\epsilon = 0.01*h$. For each simulation we let $\rho = \sup \set{\norm{(x_1(t), x_2(t))}_{\infty} \colon t \in \sp{t_{\infty}, T}}$, and use this metric to measure the convergence of our relaxed trajectories to the Zeno point. We plot the results in Figure~\ref{fig:dp} f). Note that we do not provide theoretical guarantees of the rate of convergence for this example, since the vector field is parallel to the transition surface at the origin (the Zeno accumulation point). However, the plot in Figure~\ref{fig:dp} f) nevertheless demonstrates that we have (near) linear convergence as we take $h$ and $\epsilon$ to zero, under the $\rho$ metric. We are currently working to provide formal guarantees for the rate of convergence in such cases. 

\section{conclusion}
In this paper we developed a novel class of relaxations, which we used to construct unique, well defined solutions for hybrid systems, even past the point of Zeno. The trajectories of our hybrid systems were shown to be Lipschitz continuous with respect to initial conditions and inputs, and naturally gave rise to a broad class of provably convergent discretization schemes. We provided several numerical examples, wherein we were able to accurately simulate and performed sensitivity analysis on Zeno executions. While further work is needed to extend our current framework, including the addition of non-linear guards and reset maps as well as overlapping guards, it is our conviction that the techniques developed here will provide an avenue to extract further important systems theoretic properties from hybrid systems. Moreover, we are currently working to extend our sensitivity analysis techniques to trajectories undergoing multiple transitions, with the intention of using these techniques to assess the stability of periodic orbits in hybrid systems. Such an approach has many practical applications, including finding stable periodic gates for dynamic walking robots \cite{westervelt2003hybrid}.

\bibliographystyle{acm}
\bibliography{paper}
\clearpage
\appendix
 \section{Proofs}\label{sec:proofs}
 For each of the following proofs, we provide the main arguments, omitting some details in the interest of brevity. 
 
 \subsection{Proof of Theorem \ref{thm:regular_convergence}}\label{proof:regular_convergence}
We supply the proof for the case where $|\mathcal{N}_j| = 1$; the generalization to the case where $D_j$ has multiple, disjoint guard sets is straightforward. We first demonstrate that  the claim holds for all input functions $\hat{u} \in PCD \p{\sp{0,T} , U}$, where $PCD$ denotes the class of piecewise continuously differentiable functions. We transform the system $\deriv{}{t}\bar{x} = \bar{f}_j(\bar{x},\hat{u})$ into the autonomous system defined by $(\dot {\tilde{x}}, \dot z) = \tilde{f}_j(\tilde{x},z) = (\bar{f}_j(\tilde{x},\hat{u}(z)), 1)$ which we endow with initial condition $(x_0, 0)$. Note that $z(t) =t, \forall t \in \sp{0,T}$, and thus $\tilde{x}(t) =\bar{x}(t)$, $\forall t \in \sp{0,T}$. Let $\tilde{f}_j^\epsilon$ be the $\epsilon$-relaxation of $\tilde{f}_j$, and let $(\tilde{x}^\epsilon, z)$ be the resulting trajectory, with initial data $(x^\epsilon, 0)$. Next, note that $\hat{u}$ must be non differentiable on a finite number of points $0 =t_1 < t_2 < ... <t_p = T$, $p \in \N$. Thus, on each interval $(t_i, t_{i+1}), \forall i =1,2, \dots, p-1$, $\tilde{f}_j$ is continuously differentiable in $z$. Thus, restricting both trajectories to the time interval $\sp{t_1, t_2}$, we have $\norm{(\tilde{x},z)|_{\sp{t_1,t_2}}  - (\tilde{x}^\epsilon,z)|_{\sp{t_1,t_2}}}_{\infty} = O(\epsilon)$ for each $\epsilon < \epsilon_0$ for some $\epsilon_0 >0$, by an argument similar to Lemma 2 of \cite{fiore2016contraction}. Thus, by a straight forward inductive argument we obtain $\norm{(\tilde{x},z)- (\tilde{x}^\epsilon,z)}_{\infty} = O(\epsilon)$, and thus $\norm{\bar{x} - x^\epsilon}_\infty \leq C\epsilon$, for some $C>0$. The result for our desired $u \in \BVU$ follows from noting that $PCD\p{\sp{0,T}, U}$ is dense in $\BVU$ under the $L^2$ norm, and thus we may choose $\hat{u}$ to be arbitrarily close to the desired input $u$.

\subsection{Proof of theorem \ref{thm:irregular_convergence}} First, note that $f_j^\epsilon$ is continuously differentiable in $\epsilon$ for each $\epsilon >0$, since it is constructed using a finite number of compositions and multiplications of functions which are each continuously differentiable in $\epsilon$. Thus, $\pderiv{f_j^\epsilon}{\epsilon}$ must be Lipschitz continuous for each $\epsilon \in \sp{\ubar{\epsilon}, \bar{\epsilon}}$, where $\bar{\epsilon} > \ubar{\epsilon}> 0$, since continuous functions are Lipschitz on compact domains. Thus, by Lemma 5.6.7 of \cite{polak2012optimization}, $x^\epsilon(t)$ is a Lipschitz continuous function of $\epsilon$, $\forall t \in \sp{0,T}$ and $\epsilon \in (\ubar{\epsilon} , \bar{\epsilon})$. Thus, as $\epsilon \downarrow \ubar{\epsilon}$, $x^\epsilon$ must converge uniformly to some uniformly continuous function $x^{\ubar{\epsilon}} \colon \sp{0,T} \to \bar{D}_j^{\ubar{\epsilon}}$. The desired result follows by noting that $\ubar{\epsilon}$ may be chosen to be arbitrarily small.

\subsection{Proof of Theorem \ref{thm:convergence_1}}In this case when no transitions occur, the result follows from the uniqueness of trajectories on our continuous domains. Suppose now that $x$ undergoes one transition along edge $e = (j,j')$. We may represent the trajectories of both $x$ and $x^\epsilon$ through this transition using the domain $\bar{D}_j^\epsilon$. In particular, let $\gamma$ be the solution to $\deriv{}{t}\gamma = \bar{f}_j(\gamma,u)$ with initial condition $\gamma(0)= x_0$, and let $\gamma^\epsilon$ be defined by $\deriv{}{t} \gamma^\epsilon = f_j^\epsilon(\gamma^\epsilon, u)$ with initial condition $\gamma^\epsilon(0) = x_0$. 
By Theorem~\ref{thm:regular_convergence}, $\norm{\gamma - \gamma^\epsilon}_{\infty} \leq C \epsilon$, where $\epsilon_0 ,C>0$, and $\epsilon < \epsilon_0$. For all $t$ such that $\gamma(t), \gamma^\epsilon(t) \in D_j$, we immediately have that $\indmetric{\Me}(x(t), x^\epsilon(t)) \leq C \epsilon$. For all $t$ such that $\gamma(t) \in D_e$ and $\gamma^\epsilon(t) \in D_e^\epsilon$ (i.e. when $x$ and $x^\epsilon$ have both transitioned to mode $j'$), bound $\indmetric{\Me}(x(t), x^\epsilon(t)) = \norm{x(t) -  x^\epsilon(t) } = \norm{\bar{R}_e(\gamma(t)) - \bar{R}_e^\epsilon(\gamma^\epsilon(t)) }$ $= \norm{\bar{A}_e(\gamma(t) - \gamma^\epsilon(t)) - \epsilon \hat{g}_{e'} } \leq \bar{C}\epsilon$, for some $\bar{C} >0$. If $\gamma(t) \in D_e$ but $\gamma^\epsilon(t) \in D_j^\epsilon$ (so that $x(t) \in D_{j'}$ but $x^\epsilon(t) \in D_{j}^\epsilon$ ), then by an application of the triangle inequality on $\indmetric{\Me}$, we may bound $\indmetric{\Me}(x(t), x^\epsilon(t)) \leq \norm{\gamma(t) - \gamma^\epsilon(t)} + \norm{\bar{R}_e(\gamma(t)) - \bar{R}_e^\epsilon(\gamma^\epsilon(t))} \leq \tilde{C} \epsilon$, for some $\tilde{C} >0$. The case where $\gamma(t) \in D_j$ but $\gamma^\epsilon(t)\in D_e^\epsilon$ follows similarly. Finally it is important to note that, while $x$ was transitioning, $x^\epsilon$ may have transitioned back and forth along $e$ and its partner $e'$ multiple times, yet, as a consequence of Lemma~\ref{lemma:equiv3}, $\gamma^\epsilon$ nevertheless fully captures the behavior of $x^\epsilon$ near $S_e^\epsilon \slash S_{e'}^\epsilon$. For trajectories where $x$ undergoes multiple transitions, Theorem~\ref{thm:regular_convergence} may be called inductively to complete the proof.

\section{Non-reversible Edges}\label{sec:non_reversible}
In this section we demonstrate how to extend our framework to encompass non-reversible edges. In the interest of brevity, we show how this may be done for a unimodal hybrid systems with one continuous domain $D$ (and vector field $f$), and a single non-reversible edge $e$. However, the generalization to more complicated hybrid systems with non-reversible edges is straightforward. In order to avoid needlessly introducing a large amount of slightly modified notation, we only outline these additional techniques. In particular, we demonstrate how to construct relaxed executions for these hybrid systems, and then discuss when and how the convergence theorems from the main document apply here, but do not explicitly construct the corresponding switched systems, as their structure will become apparent from the relaxations. 

We proceed by noting that $R_e(G_e)\subset \partial D$ by Definition 3.1. Since $D$ is a convex polytope, there exists a unit vector $\hat{h}_e \in \R^n$ and scalar $d_e$ such that 
\begin{equation}
R_e(G_e) \subset \set{x \in \partial D \colon h_e(x) = \hat{h}_e^Tx - d_e = 0},
\end{equation}
where by convention we choose $\hat{h}_e$ such that it points out of $D$ along $R_e(G_e)$ -- that is, $h(x) \leq 0$, $\forall x \in D$. Note, we do not assume that $R_e(G_e) \cap G_e = \emptyset$. We now define the map $\bar{R}_e^\epsilon \colon \R^n \to \R^n$ by
\begin{equation}
\bar{R}^\epsilon_e(x) = R_e(P_e(x)) - \hat{h}_eg^\epsilon_e(x),
\end{equation}
simply replacing $\hat{g}_{e'}$ with $\hat{h}_e$ when defining $\bar{R}_e$. In this case, we may now simplify
\begin{equation}
\bar{R}_e^\epsilon(x) = \bar{A}_ex + \bar{b}_e^\epsilon,
\end{equation}
where we now have that $\bar{A}_e = A_e(I - \hat{g}_e \hat{g}_e^T) - \hat{h}_e\hat{g}_{e}^T$ and $\bar{b}_e^\epsilon= A_e\hat{g}_ec_e +b_e +  \hat{h}_e(c_e+\epsilon)$. If it is the case that $\bar{A}_e$ is full rank, then we may construct relaxed transitions along $e$ using the same procedure as in the main document. That is, we use the vector fields $f_e^\epsilon$, and map $\bar{R}_e^\epsilon$ to construct relaxed transitions.

However, when $\bar{A}_e$ is not full rank these objects are ill-defined. Specifically, when $\bar{A}_e$ is not full rank we cannot use $\bar{A}_e$ to project the component of $f$ lying in the subspace $range(\bar{A}_e)^{\perp}$ back through $e$. \footnote{For example, the double pendulum with a mechanical stop from Section~\ref{sec:examples} falls into this category when $c = 0$; in particular, for this case we cannot project the differential equation for $\dot{\theta}_2$ back through the edge $\bar{e}$ of this hybrid system as in this case (if we arrange the state $x = \p{\theta_1, \dot{\theta}_1, \theta_2, \dot{\theta}_2}$) then $\hat{g}_e = \hat{h}_e = \p{0,0,-1,0}^T$ and the matrix $\bar{A}_{\bar{e}}$ works out to be

\begin{equation}
\bar{A}_{\bar{e}} = \smats{1 & 0 & 0 & 0 \\ 0 & 1 & 0 & k \\ 0 & 0 & -1 & 0 \\ 0 & 0 & 0 & 0}.
\end{equation}
}
In order to overcome this deficiency, let $\set{v_e^i}_{i =1}^{p_e}$ be a basis for $range(\bar{A}_e)^{\perp}$, where $p_e = n - rank(\bar{A}_e)$ \footnote{In the case of the double pendulum where $c= 0$, we have that $p_{\bar{e}} = 1$ and we choose $v_{\bar{e}}^1 = \p{0,0,0,1}^T$.}. In order to capture the flow of $f$ along $span\set{\set{v_e^i}_{i = 1}^{p_e}}$, we add the auxiliary state $z \in \R^{p_e}$ to our continuous state space when in mode $j$, and now define $\tilde{R}^\epsilon_e \colon \R^n \times \R^{p_e} \to \R^n$ by

\begin{equation}
\tilde{R}_e(x,z) = \bar{R}_e^\epsilon(x) + \smats{v_e^1 | \dots | v_e^{p_e}}z,
\end{equation}
which we may reformulate into
\begin{equation}
\tilde{R}_e(x,z) = \sp{ \bar{A}_e | v_e^1 | \dots | v_e^{p_e}} \smats{x \\ z} + \bar{b}_e^\epsilon,
\end{equation}
and then define $\tilde{A}_e = \sp{ \bar{A}_e | v_e^1 | \dots | v_e^{p_e}} \in \R^{n \times(n+p_e)}$, which is surjective by construction, since $range(\tilde{A}_e) = \R^n$. \footnote{For the double pendulum when $c =0$ we arrive at \begin{equation}\tilde{A}_{\bar{e}} =  \smats{1 & 0 & 0 & 0 & 0\\ 0 & 1 & 0 & k(1+c) &0 \\ 0 & 0 & -1 & 0  &0\\ 0 & 0 & 0 & 0 & 1}
\end{equation}, which is full rank and can thus be used to project the differential equation for $\dot{\theta}_2$ back through $\bar{e}$} We will now employ the right inverse of, $\tilde{A}_e$, namely $\tilde{A}_e^{\dagger} \in \R^{(p_e + n) \times n}$,  to project the dynamics of $f$ back through edge $e$, and capture this flow during a relaxed transition using the augmented state $(x^\epsilon,z) \in \R^{n+p_e}$. For the rest of the section, let $\vec{0}$ denote the $p_e$-dimensional zero vector. When we begin a relaxed execution, we will instantiate $z = \vec{0}$, and we will reset $z \to \vec{0}$ whenever a relaxed transition occurs, for reasons that will become clear momentarily. We can now define the analogue to $D_e^\epsilon$, 

\begin{equation}
\hat{D}_e^\epsilon = \set{(x,z) \in \R^n\times \R^{p_e} \colon \tilde{R}_e^\epsilon(x,z) \in D }.
\end{equation}
Next, we define
\begin{equation}
M = \sup \set{\norm{z}_{\infty} \colon \exists x \in \R^n \text{ such that } \tilde{R}_e^\epsilon(x,z) \in D},
\end{equation}
and then define the analogue to $D^\epsilon$,

\begin{equation}
\hat{D}^\epsilon = D^\epsilon \times \sp{-M, M}^{p_e},
\end{equation}
and finally the analogue to $S_e^\epsilon$
\begin{equation}
\hat{S}_e^\epsilon = S_e^\epsilon \times \sp{-M ,M}^{p_e}.
\end{equation}
That is, we confine the auxiliary state to $z$ to $\set{-M , M}^{p_e}$, so that our augmented continuous domain remains compact, but we allow $z$ to be large enough such that we can capture the full scope of $\hat{D}_e^\epsilon$ using this extra variable. Next, we define the augmented guard set 
\begin{equation}
\hat{G}_e^\epsilon = \set{(x,z) \in \tilde{S}_e^\epsilon \colon x \in G_e^\epsilon},
\end{equation}
which will triggers a discrete transition when crossed and the state is reset according to $\hat{R}_e^\epsilon \colon \R^{n + p_e} \to \R^{n+p_e}$, \footnote{Note we have overridden the original definition of $\hat{R}_e^\epsilon$ from the main document. } where 
\begin{equation}
\hat{R}_e^\epsilon (x,z) = \smats{\tilde{R}_e^\epsilon(x,z) \\ \vec{0}}.
\end{equation}

Note, a discrete transition occurs when $x \in G_e^\epsilon$, and does not depend on the value of $z$. Moreover, after the transition occurs, $z$ is reset $\vec{0}$ so that it is ready to simulate the next transition along $e$. 

Finally we are ready to define the relaxed vector field $\hat{f}^\epsilon \colon \hat{D}^\epsilon \cup \tilde{D}_e^\epsilon \times U \to \R^{n + p_e}$ by

\begin{multline}\label{eq:f_tilde}
\hat{f}^\epsilon((x,z), u) = (1 - \phi_e^\epsilon(x))\smats{f(x,u) \\ \vec{0} } \\ + \phi_e^\epsilon(x)\tilde{A}_e^{\dagger}f(\tilde{R}_e^\epsilon(x,z),u),
\end{multline}
which may be shown to be continuously differentiable. Note that, under this vector field, when $x \in D$, $\phi_e^\epsilon(x) = 0$ and $\deriv{}{t} (x, z) =(f(x,u), \vec{0})$, thus the auxiliary $z$ state does not affect the evolution of the original state $x$ when $(x,z) \notin \tilde{S}_e^\epsilon$.
That is, the auxiliary state $z$ remains dormant until the real state reaches $S_e^\epsilon$, and then $z$ begins to flow, capturing the component of $f$ that lies in $range(\bar{A}_e)^{\perp}$, as $x$ traverses $S_e^\epsilon$. Finally note that whenever $(x, z) \in \tilde{D}_e^\epsilon$ (and $g_e(x) \geq \epsilon$), the vector field $\hat{f}^\epsilon$ returns $\tilde{A}_e^{\dagger}f(\tilde{R}_e^\epsilon(x,z),u)$, which leads to the following result. 
 
 \begin{lemma}
 Let $\hat{f}^\epsilon$ be defined as in~\eqref{eq:f_tilde}. Then $\forall (x,z) \in \tilde{D}_e^\epsilon$, if we take $\deriv{}{t}(x,z) = \hat{f}^\epsilon((x,z), u)$ then we have that $\deriv{}{t} \hat{R}_e^\epsilon(x,z) = (f(\tilde{R}_e^\epsilon(x,z), u) , \vec{0})$.
 \end{lemma}
To prove the claim we compute
\begin{align}
\deriv{}{t}\hat{R}_e^\epsilon(x,z) &= \smats{\tilde{A}_e\hat{f}^\epsilon((x,z),u)\\ \vec{0} }\\
& = \smats{\tilde{A}_e \tilde{A}_e^{\dagger}f(\tilde{R}_e^\epsilon(x,z), u) \\ \vec{0} }\\
&= \smats{f(\tilde{R}_e^\epsilon(x,z), u)  \\ \vec{0}}.
\end{align}

Consequently, we can use the vector field $\hat{f}^\epsilon$ and the map $\hat{R}_e^\epsilon$ to keep track of how $(x^\epsilon,z)$ evolves during a relaxed transition. Note, that even though we are not projecting the dynamics for $z$ back through the edge $e$, since it is always reset to a value of zero and has trivial dynamics when in $\hat{D}^\epsilon$, this is not needed to keep track of how $z$ will evolve immediately after the transition.

Concretely, if the real state is instantiated at $x_0 \in D$, in order to describe a relaxed transition along $e$,  we simulate the auxiliary curve $\hat{\gamma}^\epsilon \colon \sp{0, T} \to \hat{D}^\epsilon \cup \tilde{D}_e^\epsilon$ defined by $\deriv{}{t}\hat{\gamma}^\epsilon = \hat{f}^\epsilon(\hat{\gamma}^\epsilon, u)$ with initial condition $(x_0, \vec{0})$, allowing the curve to flow into $\tilde{D}^\epsilon_e$. We then interpret $(x^\epsilon(t), z(t)) = \hat{\gamma}^\epsilon(t)$, $\forall t$ such that $\hat{\gamma}^\epsilon(t) \in \hat{D}^\epsilon$, and interpret $(x^\epsilon(t) ,z(t))  = \hat{R}_e^\epsilon(\hat{\gamma}^\epsilon(t)) = (\tilde{R}_e^\epsilon(x,z), \vec{0})$, $\forall t$ such that $\hat{\gamma}^\epsilon(t) \in \tilde{D}_e^\epsilon$. 

It is straightforward to show that analogues to Theorems~\ref{thm:regular_convergence} and \ref{thm:irregular_convergence} hold when studying the convergence of the trajectories of $\hat{f}^\epsilon$ as we take $\epsilon \to 0$; that is, sliding solutions arise when applicable and the trajectories of $\hat{f}^\epsilon$ always converge to a unique well defined limit as we take $\epsilon \to 0$. In some cases, such as the double pendulum with mechanical stop, it makes physical sense for the trajectory to "get stuck" in the relaxed strip for some time. However, in general we leave it to the practitioner to interpret this behavior. 

In order to discuss analogues to Theorems~\ref{thm:convergence_1} and \ref{thm:convergence_2}, we must first settle on the topology for the class of hybrid systems we consider here. In particular, we now define our relaxed hybrid quotient space by
\begin{equation}
\Me = \frac{D}{\Lambda_{\hat{R}_e^\epsilon}}.
\end{equation}
In order to construct trajectories over $\Me$ with multiple transitions, we further modify the construction from \cite{burden2015metrization}.

\begin{definition}\label{def:relaxed_execution2}
An execution for a relaxed unimodal hybrid dynamical system with a single non-reversible edge, given data $x_0 \in D$ and $u \in \BVU$, denoted $(x^\epsilon,z) \colon \sp{0,T} \to \M^\epsilon$ is constructed via the following algorithm.
\begin{enumerate}
\item Set $x(0)= x_0$, $z(0) =0$ and $t=0$.
\item Simulate the differential equation $\dot{\hat{\gamma}}^\epsilon(s) = \hat{f}^\epsilon(\hat{\gamma}^\epsilon(s),u(s))$ forward in time with initial condition $\hat{\gamma}^\epsilon(t) = (x^\epsilon(t), z(t))$ until time $t' = \min\set{T, \inf\set{s \colon \hat{\gamma}^\epsilon(s^+) \notin  \hat{D}^\epsilon}}$.
\item If $t' = T$ or $\hat{\gamma}^\epsilon(t') \notin \hat{G}_e^\epsilon$, let $(x^\epsilon(s) , z(s))= \hat{\gamma}^\epsilon(s), \ \forall s \in \sp{t,t'}$. Then terminate the execution. 
\item Otherwise we have that $\hat{\gamma}^\epsilon(s) \in \hat{G}_e^\epsilon$. For each $s \in [t ,t')$ set $(x^\epsilon(s), z(s)) = \hat{\gamma}^\epsilon(s)$. Set $(x^\epsilon(t'), z(t') )= \hat{R}_e^\epsilon(\hat{\gamma}^\epsilon(t'))$, set $t =t'$. Go to step 2.
\end{enumerate}
\end{definition}

Using this construction, it is straightforward to show that analogues to the results of Proposition 1, and Theorems~\ref{thm:convergence_1} and \ref{thm:convergence_2} hold for the class of hybrid systems with non-reversible edges we have considered so far in this appendix. We omit the details in the interest of brevity.

\end{document}